\documentclass[11pt]{amsart}

\usepackage{graphicx}
\usepackage[top=26mm, bottom=26mm, left=23mm, right=25mm]{geometry}
\usepackage[dvips]{color}
\usepackage{bm}
\newcommand{\vct}[1]{\bm{\mathsf{#1}}}
\newcommand{\mtx}[1]{\bm{\mathsf{#1}}}

\newcommand{\pvct}[1]{\bm{#1}}

\numberwithin{equation}{section}
\theoremstyle{definition}
\newtheorem{remark}{Remark}
\numberwithin{remark}{section}

\numberwithin{definition}{section}

\newcommand{\lsp}{\vspace{3mm}}

\renewcommand{\AA}{\mtx{A}}
\newcommand{\LL}{\mtx{L}}

\newcommand{\UU}{\mtx{U}}
\newcommand{\VV}{\mtx{V}}
\newcommand{\WW}{\mtx{W}}

\newcommand{\uu}{\vct{u}}
\newcommand{\vv}{\vct{v}}
\newcommand{\ww}{\vct{w}}

\newcommand{\pxx}{\pvct{x}}
\newcommand{\pyy}{\pvct{y}}

\begin{document}

\begin{center}
\textbf{A composite spectral scheme for variable coefficient Helmholtz problems}

\vspace{2mm}

\textit{P.G. Martinsson, Department of Applied Mathematics, University of Colorado at Boulder\\ May 31, 2012}
\end{center}

\vspace{5mm}

\begin{center}
\begin{minipage}{0.88\textwidth}
\noindent\textbf{Abstract:} A discretization scheme for variable
coefficient Helmholtz problems on two-dimensional domains is
presented.  The scheme is based on high-order spectral approximations
and is designed for problems with smooth solutions. The resulting
system of linear equations is solved using a direct solver with
$O(N^{1.5})$ complexity for the pre-computation and $O(N \log N)$
complexity for the solve.  The fact that the solver is direct is a
principal feature of the scheme, since iterative methods tend to
struggle with the Helmholtz equation.  Numerical examples
demonstrate that the scheme is fast and highly accurate.  For
instance, using a discretization with 12 points per wave-length, a
Helmholtz problem on a domain of size $100 \times 100$ wavelengths was
solved to ten correct digits. The computation was executed on an
office desktop; it involved 1.6M degrees of freedom and required 100
seconds for the pre-computation, and 0.3 seconds for the actual solve.
%
\end{minipage}
\end{center}

\section{Introduction}
\label{sec:intro}

The paper describes a technique for constructing an approximate solution to the
variable coefficient Helmholtz equation
\begin{equation}
\label{eq:basic}
\left\{\begin{aligned}
-\Delta u(\pxx) - \kappa^{2}\,(1-b(\pxx))\,u(\pxx) =&\ 0\qquad &\pxx \in \Omega,\\
u(\pxx) =&\ f(\pxx)\qquad &\pxx \in \Gamma,
\end{aligned}\right.
\end{equation}
where $\Omega$ is a rectangular domain with boundary $\Gamma$,
where the Helmholtz parameter $\kappa$ is real,
and where $b$ is a given smooth scattering potential.
The scheme can straight-forwardly
be adapted to handle other variable coefficient elliptic problems, as
well as free space scattering problems in $\mathbb{R}^{2}$. The
primary limitation of the method is that it requires the solution $u$
to be smooth in $\Omega$.

\begin{figure}
\includegraphics[height=70mm]{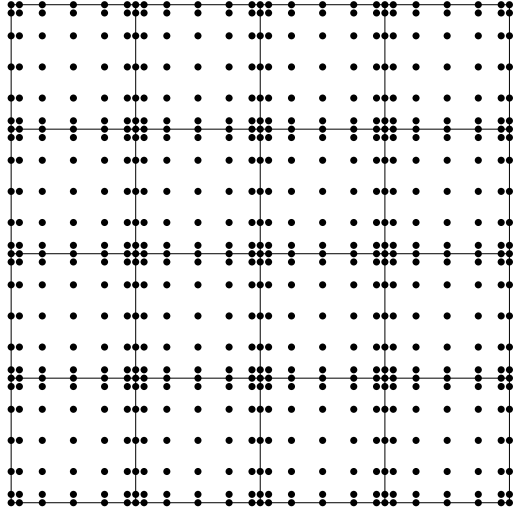}
\caption{The box $\Omega = [0,1]^{2}$ is split into $4 \times 4$ leaf boxes,
and a Cartesian grid of Chebyshev nodes is placed on each leaf box.
The figure shows local grids of size $7 \times 7$ for clarity; in actual
computations, local grids of size $21 \times 21$ were typically used.}
\label{fig:grid}
\end{figure}

The equation (\ref{eq:basic}) is discretized via a composite spectral
scheme.  The domain $\Omega$ is split into small square (or
rectangular) patches. On each patch, the solution $u$ is represented
via tabulation on a tensor product grid of Chebyshev points, see
Figure \ref{fig:grid}.  The Laplace operator is approximated via a
spectral differentiation matrix acting on each local grid, and then
equation (\ref{eq:basic}) is enforced strongly at all tabulation nodes
in the interior of each patch.  To glue patches together, continuity of
both the potential $u$ and its normal derivative are enforced at the
spectral interpolation nodes on the boundaries between patches.

The discretization scheme is combined with a direct solver for the resulting
linear system. The fact that the solver is direct rather than iterative is a
principal feature of the scheme, since iterative solvers tend to struggle for
Helmholtz problems of the kind considered here \cite{2008_duan_rokhlin}.
The direct solver organizes the patches in the discretization into a binary
tree of successively larger patches. The solver then involves two stages, one
that involves an upwards pass, and one that involves a downwards pass:
\begin{enumerate}
\item A \textit{pre-computation stage} where an approximation to the
  solution operator for (\ref{eq:basic}) is computed. This is done via
  a single sweep of the hierarchical tree, going from smaller patches
  to larger. For each leaf in the tree, a local solution operator, and
  an approximation to the Dirichlet-to-Neumann (DtN) map for the patch
  are constructed. For a parent node in the tree, a local solution
  operator and a local DtN operator are computed from an equilibrium
  equation formed using the DtN operators of the children of the
  patch.  The pre-computation stage has asymptotic complexity
  $O(N^{1.5})$.
\item A \textit{solve} stage that takes as input a vector of Dirichlet
  data tabulated on $\Gamma$, and constructs tabulated values of $u$
  at all internal grid points. The solve stage involves a single
  downwards sweep through the hierarchical tree of patches, going
  from larger patches to smaller. The solve stage has asymptotic
  complexity $O(N\log N)$.
\end{enumerate}

Numerical experiments indicate that the spectral convergence of the
method makes it both highly accurate and computationally
efficient. For instance, the equation (\ref{eq:basic}) was solved for
a box whose size exceeded $100 \times 100$ wave-lengths in less than 2
minutes on a standard office laptop.  A $20$-th order spectral scheme
with 12 points per wave-length was used in the local approximation on
the patches.  The resulting solution was accurate to between 7 and 10
digits, depending on the nature of the scattering potential $b$ in
(\ref{eq:basic}).  The discretization used a total of $N=1.6\cdot
10^{6}$ degrees of freedom.  The computational time was dominated by
the pre-computation stage; the actual solve stage took only $0.3$
seconds. This makes the scheme particularly powerful in situations
where an equation such as (\ref{eq:basic}) needs to be solved for a
sequence of different boundary functions $f$.

The scheme proposed is conceptually related to a direct solver for the
Lippman-Schwinger equation proposed in 2002 by Yu Chen
\cite{2002_chen_direct_lippman_schwinger}. The schemes are different
in that the method proposed here is not based on a Lippman-Schwinger
formulation, and uses spectral approximations on the smallest patches
in the hierarchical tree.  Comparing the efficiencies of the two
schemes is difficult since the paper
\cite{2002_chen_direct_lippman_schwinger} does not report numerical
results, and we have been unable to find reports of implementations of
the scheme.  The scheme proposed here is also conceptually related to
the classical nested dissection algorithm for finite element and
finite difference matrices \cite{george_1973}, and to recently
proposed $O(N^{1.5})$ direct solvers for BIEs on surfaces in 3D
\cite{2011_gillman_dissertation}.

For clarity, the current paper focusses on the simple boundary value
problem (\ref{eq:basic}) involving the Helmholtz elliptic operator and
Dirichlet boundary data. The scheme can with trivial modifications be
applied to more general elliptic operators
\begin{multline*}
-c_{11}(\pxx)[\partial_{1}^{2}u](\pxx)
-2c_{12}(\pxx)[\partial_{1}\partial_{2}u](\pxx)
-c_{22}(\pxx)[\partial_{2}^{2}u](\pxx)\\
+c_{1}(\pxx)[\partial_{1}u](\pxx)
+c_{2}(\pxx)[\partial_{2}u](\pxx)
+c(\pxx)\,u(\pxx) = 0,
\end{multline*}
coupled with Dirichlet, Neumann, or mixed boundary data. It has for
instance been successfully tested on convection-diffusion problems
that are strongly dominated (by a factor of $10^{4}$) by the
convection term, see Section \ref{sec:convdiff}. Moreover, the scheme
can with minor modifications be applied to a free space scattering
problem such as
\begin{equation}
\label{eq:freespace}
-\Delta u(\pxx) - \kappa^{2}\,(1-b(\pxx))\,u(\pxx) = f(\pxx), \qquad \pxx \in \mathbb{R}^{2},
\end{equation}
coupled with appropriate radiation conditions at infinity. A standard
assumption is that $f$ is supported outside of some (bounded) square
region $\Omega$ while the smooth function $b$ is supported inside
$\Omega$. The scheme described in this note computes the DtN operator
for (\ref{eq:freespace}) on $\Omega$. The DtN operator for the
exterior domain $\Omega^{\rm c}$ can be computed via Boundary Integral
Equation techniques; and by combining the two, one can solve the free
space scattering problem, see Section \ref{sec:freespace}.

The method proposed has a vulnerability in that it crucially relies on
the existence of DtN operators for all patches in the hierarchical
tree. This can be problematic due to resonances: For certain
wave-numbers $\kappa$, there exist non-trivial solutions that have
zero Dirichlet boundary data. We have found that in practice, this
problem almost never arises when processing domains that are a couple
of hundred wave-lengths or less in size. Moreover, if a resonant patch
should be encountered, this will be detected and counter-measures can
be taken, see Sections \ref{sec:resonance} and \ref{sec:reformulation}.

The asymptotic complexity of the proposed method is $O(N^{1.5})$. For the case
where the wave-number $\kappa$ is increased as $N$ grows to keep a constant number
of discretization points per wave-length (i.e.~$\kappa \sim N^{0.5}$), we do not
know how to improve the complexity. However, for the case where the wave-number is
kept constant as $N$ increases, $O(N)$ complexity can very likely be attained by
exploiting internal structure in the DtN operators. The resulting scheme would be
a spectral version of recently published accelerated nested dissection schemes such
as \cite{2007_leborne_HLU,2011_ying_nested_dissection_2D,2009_xia_superfast}.

An early version of the work reported was published on arXiv as \cite{2011_martinsson_gauss}.

The paper is organized as follows:
Section \ref{sec:prelim} introduces notation and lists some classical material on spectral
interpolation and differentiation.
Section \ref{sec:leaf} describes how to compute the solution operator and the DtN operator for a leaf in tree
(which is discretized via a single tensor-product grid of Chebyshev nodes).
Section \ref{sec:merge} describes how the DtN operator for a larger patch consisting of two small
patches can be computed if the DtN operators for the smaller patches are given.
Section \ref{sec:hierarchy} describes the full hierarchical scheme.
Section \ref{sec:num} reports the results of some numerical experiments.
Section \ref{sec:extensions} describes how the scheme can be extended to more general situations.

\section{Preliminaries --- spectral differentiation}
\label{sec:prelim}

This section introduces notation for spectral differentiation on
tensor product grids of Chebyshev nodes on the square domain
$[-a,a]^{2}$.  This material is classical, see, e.g., Trefethen
\cite{2000_trefethen_spectral_matlab}.  (While we restrict attention
to square boxes here, all techniques generalize trivially to
rectangular boxes of moderate aspect ratios.)

Let $p$ denote a positive integer. The \textit{Chebyshev nodes} on
$[-a,a]$ are the points
$$
t_{i} = a\,\cos((i-1)\pi/(p-1)),\qquad i = 1,\,2,\,3,\,\dots,\,p.
$$
Let $\{\pxx_{k}\}_{k=1}^{p^{2}}$ denote the set of
points of the form $(t_{i},\,t_{j})$ for $1 \leq i,j \leq p$.  Let
$\mathcal{P}_{p}$ denote the linear space of sums of tensor products
of polynomials of degree $p-1$ or less. $\mathcal{P}_{p}$ has
dimension $p^{2}$.  Given a vector $\uu \in \mathbb{R}^{p^{2}}$, there is a
unique function $u \in \mathcal{P}_{p}$ such that $u(\pxx_{k}) =
\uu(k)$ for $1 \leq k \leq p^{2}$. (A reason Chebyshev nodes are of interest
is that for any fixed $\pxx \in [-a,a]^{2}$, the map $\uu \mapsto
u(\pxx)$ is stable.)  Now define $\mtx{D}$, $\mtx{E}$, and $\mtx{L}$
as the unique $p^{2}\times p^{2}$ matrices such that
\begin{align}
[\mtx{D}\uu](k) =&\ [\partial_{1}u](\pxx_{k}),\qquad k = 1,\,2,\,3,\,\dots,\,p^{2},\\
[\mtx{E}\uu](k) =&\ [\partial_{2}u](\pxx_{k}),\qquad k = 1,\,2,\,3,\,\dots,\,p^{2},\\
[\mtx{L}\uu](k) =&\ [-\Delta     u](\pxx_{k}),\qquad k = 1,\,2,\,3,\,\dots,\,p^{2}.
\end{align}

\section{Leaf computation}
\label{sec:leaf}

This section describes the construction of a discrete approximation to
the Dirichlet-to-Neumann operator associated with the boundary value
problem (\ref{eq:basic}) for a square patch $\Omega$.  We discretize
(\ref{eq:basic}) via a spectral method on a tensor product grid of
Chebyshev nodes on $\Omega$. In addition to the DtN operator, we also
construct a solution operator to (\ref{eq:basic}) that maps the
Dirichlet data on the nodes on the boundary of $\Omega$ to the value
of $u$ at all internal interpolation nodes.

\subsection{Notation}
\begin{figure}
\begin{tabular}{ccc}
\includegraphics[height=50mm]{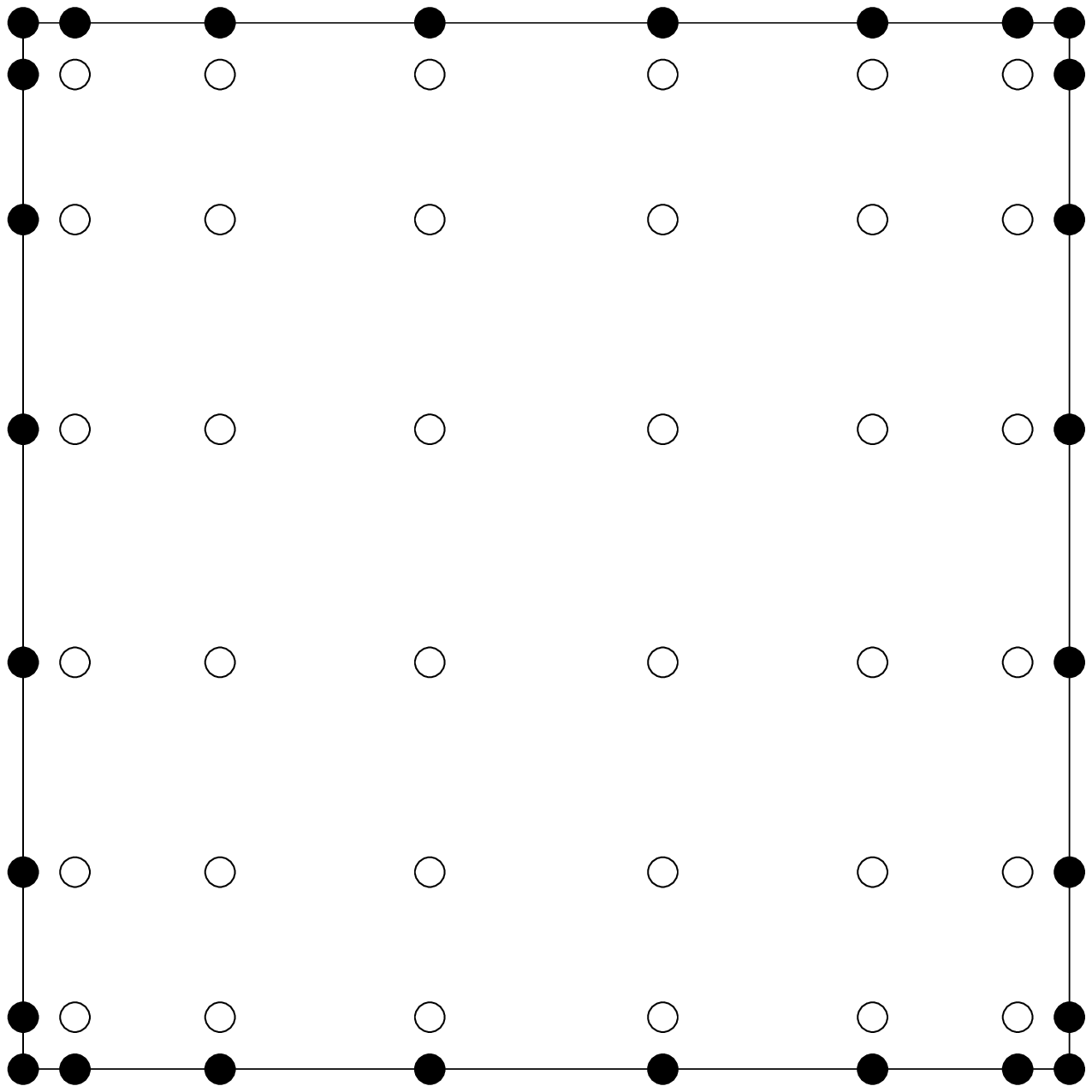}
& \mbox{}\hspace{25mm}\mbox{} &
\includegraphics[height=50mm]{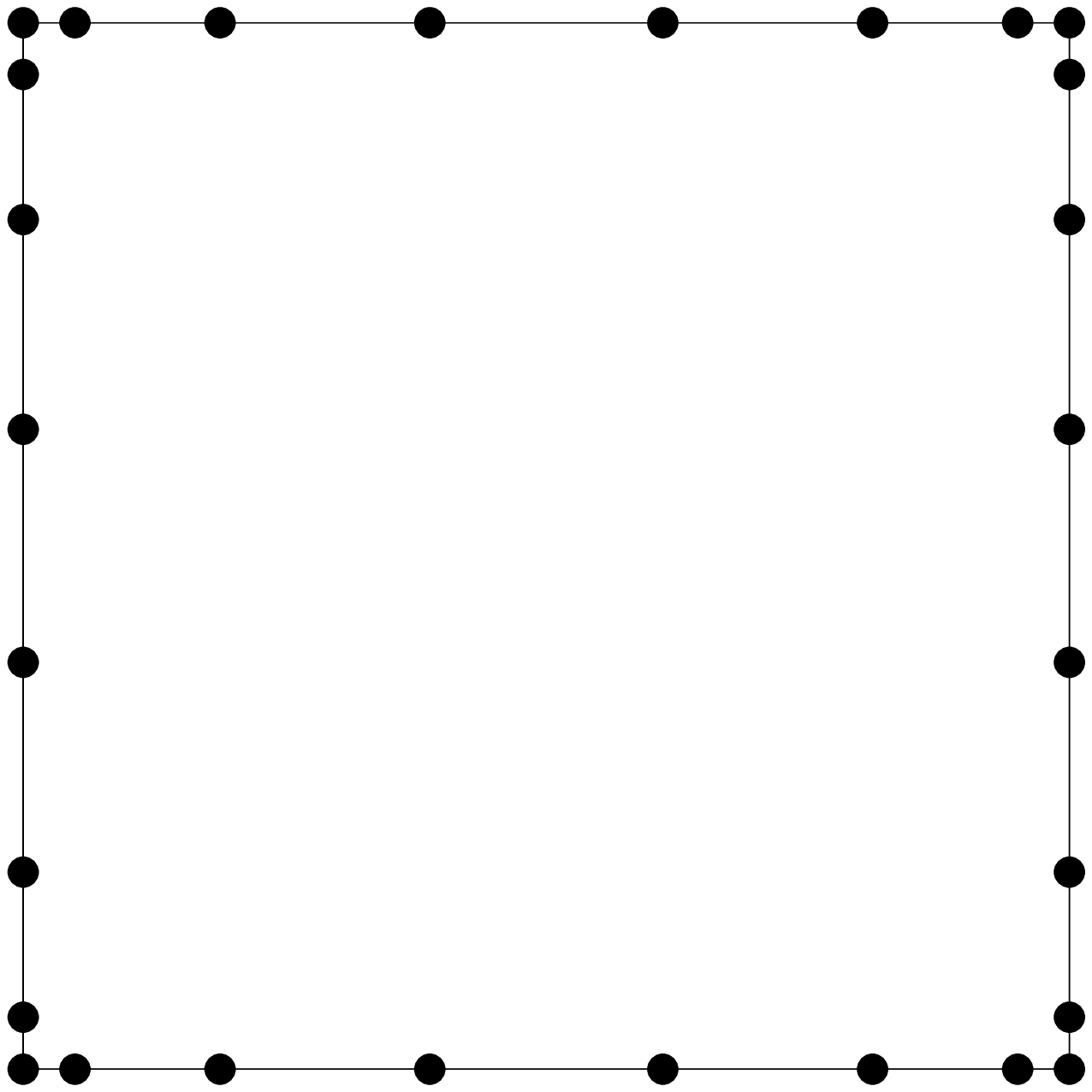}\\
(a) && (b)
\end{tabular}
\caption{Notation for the leaf computation in Section \ref{sec:leaf}.
(a) A leaf before elimination of interior (white) nodes.
(b) A leaf after elimination of interior nodes.}
\label{fig:leaf}
\end{figure}

Let $\Omega$ denote a square patch.  Let $\{\pxx_{k}\}_{k=1}^{p^{2}}$
denote the nodes in a tensor product grid of $p\times p$ Chebyshev
nodes.  Partition the index set
$$
\{1,\,2,\,\dots,\,p^{2}\} = I_{\rm e}\cup I_{\rm i}
$$
in such a way that $I_{\rm e}$ contains all nodes on the boundary of $\Omega$, and
$I_{\rm i}$ denotes the set of interior nodes, see Figure \ref{fig:leaf}(a).
Let $u$ be a function that satisfies (\ref{eq:basic}) on $\Omega$ and let
$$
\uu = [u(\pxx_{k})]_{k=1}^{p^{2}},\qquad
\vv = [\partial_{1}u(\pxx_{k})]_{k=1}^{p^{2}},\qquad
\ww = [\partial_{2}u(\pxx_{k})]_{k=1}^{p^{2}},
$$
denote the vectors of samples of $u$ and its partial derivatives.
We define the short-hands
$$
\uu_{\rm i} = \uu(I_{\rm i}),\qquad
\vv_{\rm i} = \vv(I_{\rm i}),\qquad
\ww_{\rm i} = \ww(I_{\rm i}),\qquad
\uu_{\rm e} = \uu(I_{\rm e}),\qquad
\vv_{\rm e} = \vv(I_{\rm e}),\qquad
\ww_{\rm e} = \ww(I_{\rm e}).
$$
Let $\mtx{L}$, $\mtx{D}$, and $\mtx{E}$ denote spectral differentiation matrices corresponding to
the operators $-\Delta$, $\partial_{1}$, and $\partial_{2}$, respectively (see Section \ref{sec:prelim}).
We use the short-hand
$$
\mtx{D}_{\rm i,e} = \mtx{D}(I_{\rm i},I_{\rm e})
$$
to denote the part of the differentiation matrix $\mtx{D}$ that maps exterior nodes to interior nodes,
etc.

\subsection{Equilibrium condition}
The operator (\ref{eq:basic}) is approximated via the matrix
$$
\AA = -\LL - \kappa^{2}\mbox{diag}(\vct{b}),
$$
where $\vct{b}$ denotes the vector of pointwise values of $b$,
$$
\vct{b} = [b(\pvct{x}_{k})]_{k=1}^{p^{2}}.
$$
The equation we enforce on $\Omega$ is that the vector $\AA\,\uu$ should evaluate to zero \textit{at all internal nodes},
\begin{equation}
\label{eq:A_interior}
\AA_{\rm i,i}\,\uu_{\rm i} +
\AA_{\rm i,e}\,\uu_{\rm e} = \vct{0},
\end{equation}
where
$$
\AA_{\rm i,i} = \AA(I_{\rm i},I_{\rm i}),\qquad
\AA_{\rm i,e} = \AA(I_{\rm i},I_{\rm e}).
$$
Solving (\ref{eq:A_interior}) for $\uu_{\rm i}$, we obtain
\begin{equation}
\label{eq:reconstruction_leaf}
\uu_{\rm i} = \UU\,\uu_{\rm e},
\end{equation}
where
\begin{equation}
\label{eq:U_leaf}
\UU = -\bigl(\AA_{\rm i,i}\bigr)^{-1}\,\AA_{\rm i,e}.
\end{equation}

\subsection{Constructing the DtN operator}
Let $\VV$ and $\WW$ denote the matrices that map boundary values of
the potential to boundary values of $\partial_{1}u$ and
$\partial_{2}u$.  These are constructed as follows: Given the
potential $\uu_{\rm e}$ on the boundary, we reconstruct the potential
$\uu_{\rm i}$ in the interior via
(\ref{eq:reconstruction_leaf}). Then, since the potential is known on
all Chebyshev nodes in $\Omega$, we can determine the gradient on the
boundary $\{\vv_{\rm e},\,\ww_{\rm e}\}$ via spectral differentiation
on the entire domain. To formalize, we find
$$
\vv_{\rm e} =
\mtx{D}_{\rm e,e}\,\uu_{\rm e} +
\mtx{D}_{\rm e,i}\,\uu_{\rm i}
=
\mtx{D}_{\rm e,e}\,\uu_{\rm e} +
\mtx{D}_{\rm e,i}\,\UU\,\uu_{\rm e}
=
\VV\,\uu_{\rm e},
$$
where
\begin{equation}
\label{eq:V_leaf}
\VV = \mtx{D}_{\rm e,e} + \mtx{D}_{\rm e,i}\,\UU.
\end{equation}
An analogous computation for $\ww_{\rm e}$ yields
\begin{equation}
\label{eq:W_leaf}
\WW = \mtx{E}_{\rm e,e} + \mtx{E}_{\rm e,i}\,\UU.
\end{equation}


\section{Merge operation}
\label{sec:merge}

\begin{figure}
\setlength{\unitlength}{1mm}
\begin{picture}(95,55)
\put(00,00){\includegraphics[height=55mm]{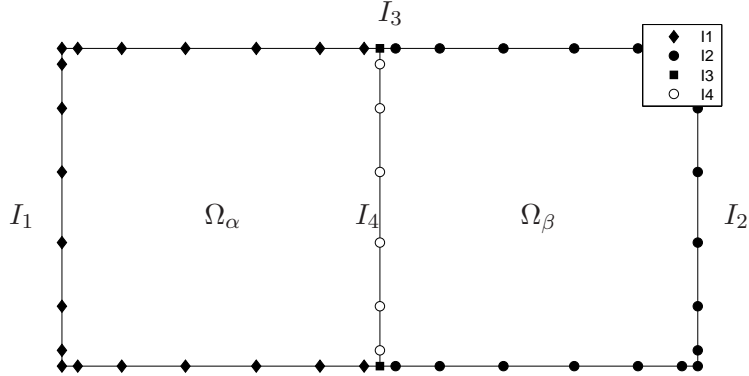}}
\put(24,25){$\Omega_{\alpha}$}
\put(66,25){$\Omega_{\beta}$}
\put(-2,25){$I_{1}$}
\put(93,25){$I_{2}$}
\put(44,25){$I_{4}$}
\put(47,52){$I_{3}$}
\end{picture}
\caption{Notation for the merge operation described in Section \ref{sec:merge}.
The rectangular domain $\Omega$ is formed by two squares $\Omega_{\alpha}$ and $\Omega_{\beta}$.
The sets $I_{1}$, $I_{2}$, and $I_{3}$ form the exterior nodes (black), while
$I_{4}$ consists of the interior nodes (white).}
\label{fig:siblings_notation}
\end{figure}

Let $\Omega$ denote a rectangular domain consisting of the
union of the two smaller rectangular domains,
$$
\Omega = \Omega_{\alpha} \cup \Omega_{\beta},
$$ as shown in Figure \ref{fig:siblings_notation}. Moreover, suppose
that approximations to the DtN operators for $\Omega_{\alpha}$ and
$\Omega_{\beta}$ are available.  (Represented as matrices that map
boundary values of $u$ to boundary values of $\partial_{1}u$ and
$\partial_{2}u$.) This section describes how to compute a solution
operator $\UU$ that maps the value of a function $u$ that is tabulated
on the boundary of $\Omega$ to the values of $u$ on interpolation
nodes on the internal boundary, as well as operators $\VV$ and $\WW$
that map boundary values of $u$ on the boundary of $\Omega$ to values
of the $\partial_{1}u$ and $\partial_{2}u$ tabulated on the boundary.

\subsection{Notation}
Let $\Omega$ denote a box with children $\Omega_{\alpha}$ and
$\Omega_{\beta}$.  For concreteness, let us assume that $\Omega_{\alpha}$ and
$\Omega_{\beta}$ share a vertical edge.  We partition the points on
$\partial\Omega_{\alpha}$ and $\partial\Omega_{\beta}$ into four sets:
\begin{tabbing}
\mbox{}\hspace{5mm}\= $I_{1}$ \hspace{4mm} \=
Boundary nodes of $\Omega_{\alpha}$ that are not boundary nodes of $\Omega_{\beta}$.\\
\> $I_{2}$ \> Boundary nodes of $\Omega_{\beta}$ that are not boundary nodes of $\Omega_{\alpha}$.\\
\> $I_{3}$ \> The two nodes that are boundary nodes of $\Omega_{\alpha}$, of $\Omega_{\beta}$, and also of the union box $\Omega$.\\
\> $I_{4}$ \> Boundary nodes of both $\Omega_{\alpha}$ and $\Omega_{\beta}$ that are \textit{not} boundary nodes of the union box
$\Omega$.
\end{tabbing}
Figure \ref{fig:siblings_notation} illustrates the definitions of the $I_{j}$'s.
Let $u$ denote a function such that
$$
-\Delta u(\pxx) - \kappa^{2}\,b(\pxx)\,u(\pxx) = 0,\qquad \pxx \in \Omega,
$$
and let $\uu_{j}$, $\vv_{j}$, $\ww_{j}$
denote the values of $u$, $\partial_{1}u$, and $\partial_{2}u$,
restricted to the nodes in the set ``$j$''. Moreover, set
\begin{equation}
\label{eq:defuiue}
\uu_{\rm i} = \uu_{4},
\qquad\mbox{and}\qquad
\uu_{\rm e} = \left[\begin{array}{c} \uu_{1} \\ \uu_{2} \\ \uu_{3}\end{array}\right].
\end{equation}
Finally, let $\VV^{\alpha}$, $\WW^{\alpha}$, $\VV^{\beta}$, $\WW^{\beta}$ denote the
operators that map potential values on the boundary to values of $\partial_{1}u$ and $\partial_{2}u$
on the boundary for the boxes $\Omega_{\alpha}$ and $\Omega_{\beta}$. We partition these matrices according
to the numbering of nodes in Figure \ref{fig:siblings_notation},
\begin{equation}
\label{eq:bittersweet1}
\left[\begin{array}{c}
\vv_{1}\\ \vv_{3}\\ \vv_{4}
\end{array}\right] =
\left[\begin{array}{ccc}
\mtx{V}_{1,1}^{\alpha} & \mtx{V}_{1,3}^{\alpha} & \mtx{V}_{1,4}^{\alpha} \\
\mtx{V}_{3,1}^{\alpha} & \mtx{V}_{3,3}^{\alpha} & \mtx{V}_{3,4}^{\alpha} \\
\mtx{V}_{4,1}^{\alpha} & \mtx{V}_{4,3}^{\alpha} & \mtx{V}_{4,4}^{\alpha}
\end{array}\right]\,
\left[\begin{array}{c}
\uu_{1}\\ \uu_{3}\\ \uu_{4}
\end{array}\right],
\qquad
\left[\begin{array}{c}
\ww_{1}\\ \ww_{3}\\ \ww_{4}
\end{array}\right] =
\left[\begin{array}{ccc}
\mtx{W}_{1,1}^{\alpha} & \mtx{W}_{1,3}^{\alpha} & \mtx{W}_{1,4}^{\alpha} \\
\mtx{W}_{3,1}^{\alpha} & \mtx{W}_{3,3}^{\alpha} & \mtx{W}_{3,4}^{\alpha} \\
\mtx{W}_{4,1}^{\alpha} & \mtx{W}_{4,3}^{\alpha} & \mtx{W}_{4,4}^{\alpha}
\end{array}\right]\,
\left[\begin{array}{c}
\uu_{1}\\ \uu_{3}\\ \uu_{4}
\end{array}\right],
\end{equation}
and
\begin{equation}
\label{eq:bittersweet2}
\left[\begin{array}{c}
\vv_{2}\\ \vv_{3}\\ \vv_{4}
\end{array}\right] =
\left[\begin{array}{ccc}
\mtx{V}_{2,2}^{\beta} & \mtx{V}_{2,3}^{\beta} & \mtx{V}_{2,4}^{\beta} \\
\mtx{V}_{3,2}^{\beta} & \mtx{V}_{3,3}^{\beta} & \mtx{V}_{3,4}^{\beta} \\
\mtx{V}_{4,2}^{\beta} & \mtx{V}_{4,3}^{\beta} & \mtx{V}_{4,4}^{\beta}
\end{array}\right]\,
\left[\begin{array}{c}
\uu_{2}\\ \uu_{3}\\ \uu_{4}
\end{array}\right],
\qquad
\left[\begin{array}{c}
\ww_{2}\\ \ww_{3}\\ \ww_{4}
\end{array}\right] =
\left[\begin{array}{ccc}
\mtx{W}_{2,2}^{\beta} & \mtx{W}_{2,3}^{\beta} & \mtx{W}_{2,4}^{\beta} \\
\mtx{W}_{3,2}^{\beta} & \mtx{W}_{3,3}^{\beta} & \mtx{W}_{3,4}^{\beta} \\
\mtx{W}_{4,2}^{\beta} & \mtx{W}_{4,3}^{\beta} & \mtx{W}_{4,4}^{\beta}
\end{array}\right]\,
\left[\begin{array}{c}
\uu_{2}\\ \uu_{3}\\ \uu_{4}
\end{array}\right].
\end{equation}

\subsection{Equilibrium condition}
\begin{figure}
\begin{tabular}{ccc}
\includegraphics[height=40mm]{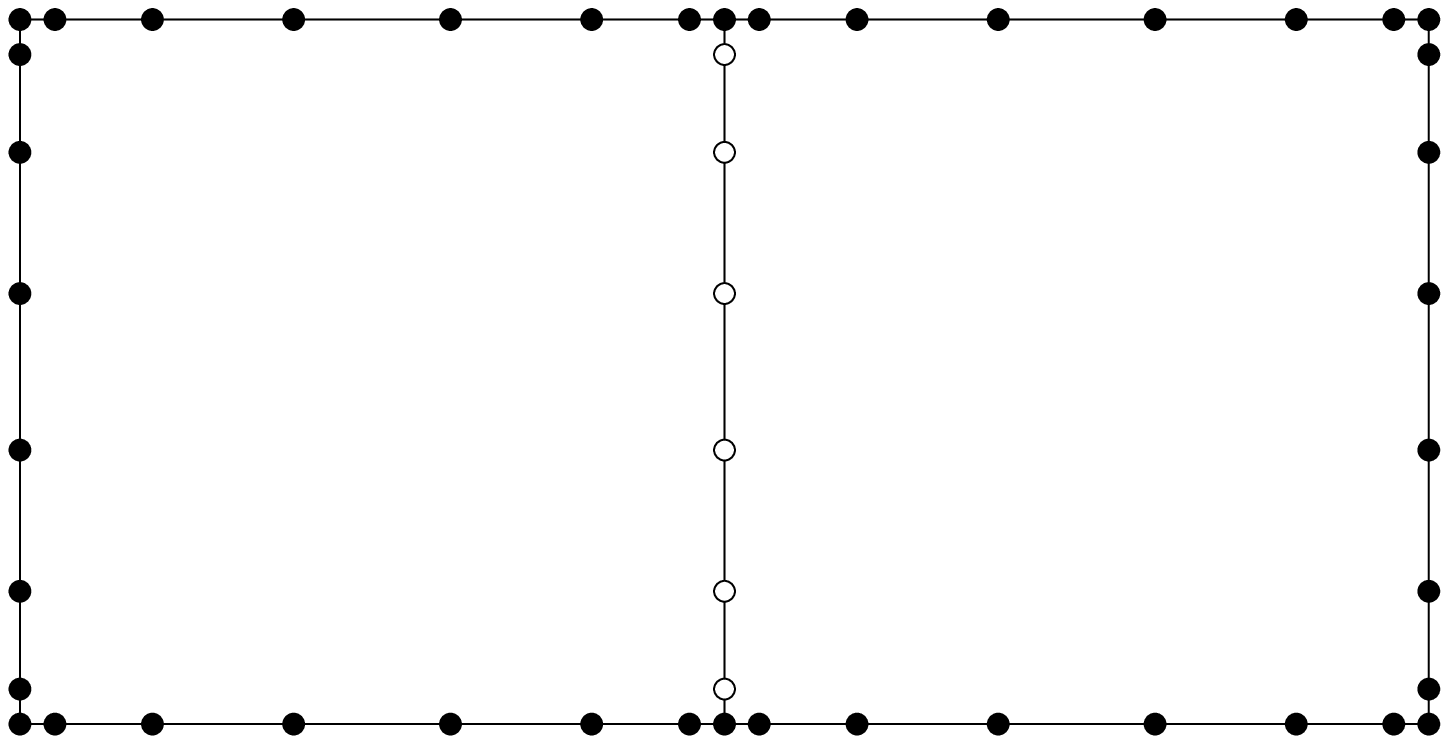}
& \mbox{}\hspace{5mm}\mbox{} &
\includegraphics[height=40mm]{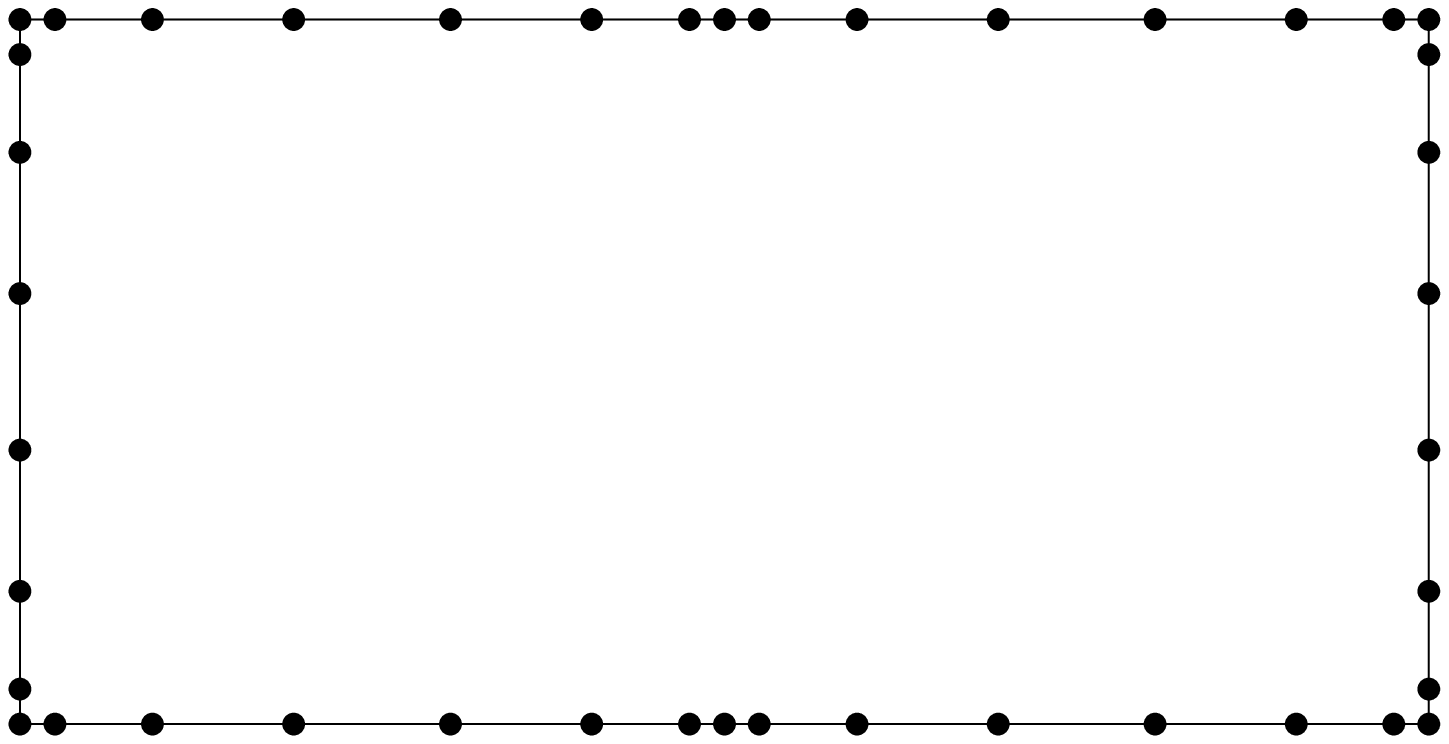}\\
(a) && (b)
\end{tabular}
\caption{Merge operation for two small boxes to form  a new large box.
(a) Before elimination of interior (white) nodes.
(b) After elimination of interior nodes.}
\label{fig:siblings}
\end{figure}

Suppose that we are given a tabulation of boundary values of a function $u$ that satisfies
(\ref{eq:basic})  on $\Omega$. In other words, we are given the vectors $\uu_{1}$, $\uu_{2}$,
and $\uu_{3}$. We can then reconstruct the values of the potential on the interior boundary
(tabulated in the vector $\uu_{4}$) by using information in (\ref{eq:bittersweet1}) and
(\ref{eq:bittersweet2}). Simply observe that there are two equations specifying the normal
derivative across the internal boundary (tabulated in $\vv_{4}$), and combine these equations:
$$
\mtx{V}_{4,1}^{\alpha}\uu_{1} + \mtx{V}_{4,3}^{\alpha}\uu_{3} + \mtx{V}_{4,4}^{\alpha}\uu_{4} =
\mtx{V}_{4,2}^{\beta}\uu_{2} + \mtx{V}_{4,3}^{\beta}\uu_{3} + \mtx{V}_{4,4}^{\beta}\uu_{4}.
$$
Solving for $\uu_{4}$ we get
\begin{equation}
\label{eq:u4}
\uu_{4} =
\bigl(\VV^{\alpha}_{4,4} - \VV^{\beta}_{4,4}\bigr)^{-1}
\bigl( \VV^{\beta }_{4,2}\uu_{2}
      +\VV^{\beta }_{4,3}\uu_{3}
      -\VV^{\alpha}_{4,1}\uu_{1}
      -\VV^{\alpha}_{4,3}\uu_{3}\bigr).
\end{equation}
Now set
\begin{equation}
\label{eq:U_parent}
\UU =
\bigl(\VV^{\alpha}_{4,4} - \VV^{\beta}_{4,4}\bigr)^{-1}
\bigl[-\VV^{\alpha}_{4,1}\ \big|\
       \VV^{\beta}_{4,2}\  \big|\
       \VV^{\beta}_{4,3} - \VV^{\alpha}_{4,3} \bigr],
\end{equation}
to find that (\ref{eq:u4}) is (in view of (\ref{eq:defuiue})) precisely the desired formula
\begin{equation}
\label{eq:desired_merge}
\uu_{\rm i} = \UU\,\uu_{\rm e}.
\end{equation}
The net effect of the merge operation is to eliminate the interior
tabulation nodes in $\Omega_{\alpha}$ and $\Omega_{\beta}$ so that only
boundary nodes in the union box $\Omega$ are kept, as illustrated in
Figure \ref{fig:siblings}.

\subsection{Constructing the DtN operators for the union box}
We will next build a matrix $\VV$ that constructs the
derivative $\partial_{1}u$ on $\partial\Omega$ given values
of $u$ on $\partial\Omega$. In other words
$$
\left[\begin{array}{c} \vv_{1} \\ \vv_{2} \\ \vv_{3} \end{array}\right] =
\VV\,
\left[\begin{array}{c} \uu_{1} \\ \uu_{2} \\ \uu_{3} \end{array}\right].
$$
To this end, observe from (\ref{eq:bittersweet1}) and (\ref{eq:bittersweet2}) that
\begin{align}
\label{eq:v1}
\vv_{1} =&\
\VV^{\alpha}_{1,1}\,\uu_{1} + \VV^{\alpha}_{1,3}\,\uu_{3} + \VV^{\alpha}_{1,4}\,\uu_{4} =
\VV^{\alpha}_{1,1}\,\uu_{1} + \VV^{\alpha}_{1,3}\,\uu_{3} + \VV^{\alpha}_{1,4}\,\UU\,\uu_{\rm e}\\
\label{eq:v2}
\vv_{2} =&\
\VV^{\beta }_{2,2}\,\uu_{2} + \VV^{\alpha}_{2,3}\,\uu_{3} + \VV^{\alpha}_{2,4}\,\uu_{4} =
\VV^{\beta }_{2,2}\,\uu_{2} + \VV^{\alpha}_{2,3}\,\uu_{3} + \VV^{\alpha}_{2,4}\,\UU\,\uu_{\rm e}.
\end{align}
Equations (\ref{eq:bittersweet1}) and (\ref{eq:bittersweet2}) provide two different formulas for
$\vv_{3}$, either of which could be used. For numerical stability, we use the average of the two:
\begin{align}
\vv_{3} =&\ \frac{1}{2}\bigl(\VV^{\alpha}_{3,1}\uu_{1} + \VV^{\alpha}_{3,3}\uu_{3} + \VV^{\alpha}_{3,4}\uu_{4} +
                             \VV^{\beta }_{3,2}\uu_{2} + \VV^{\beta }_{3,3}\uu_{3} + \VV^{\beta }_{3,4}\uu_{4}\bigr)\\
\label{eq:v3}
        =&\ \frac{1}{2}\bigl(\VV^{\alpha}_{3,1}\uu_{1} + \VV^{\alpha}_{3,3}\uu_{3} + \VV^{\alpha}_{3,4}\,\UU\,\uu_{\rm e} +
                             \VV^{\beta }_{3,2}\uu_{2} + \VV^{\beta }_{3,3}\uu_{3} + \VV^{\beta }_{3,4}\,\UU\,\uu_{\rm e}\bigr).
\end{align}
Combining (\ref{eq:v1}) -- (\ref{eq:v3}) we obtain
$$
\left[\begin{array}{c} \vv_{1} \\ \vv_{2} \\ \vv_{3} \end{array}\right] =
\left(
\left[\begin{array}{ccc}
\VV_{1,1}^{\alpha} & \mtx{0} & \VV_{1,3}^{\alpha} \\
\mtx{0} & \VV_{2,2}^{\beta } & \VV_{2,3}^{\beta }  \\
\tfrac{1}{2}\VV_{3,1}^{\alpha} & \tfrac{1}{2}\VV_{3,2}^{\beta} & \tfrac{1}{2}\VV_{3,3}^{\alpha} + \tfrac{1}{2}\VV_{3,3}^{\beta}
\end{array}\right] +
\left[\begin{array}{c}
\VV_{1,4}^{\alpha} \\
\VV_{2,4}^{\beta}  \\
\tfrac{1}{2}\VV_{3,4}^{\alpha} + \tfrac{1}{2}\VV_{3,4}^{\beta}
\end{array}\right]\,\UU\right)\,
\left[\begin{array}{c} \uu_{1} \\ \uu_{2} \\ \uu_{3} \end{array}\right].
$$
In other words,
\begin{equation}
\label{eq:V_parent}
\VV =
\left[\begin{array}{ccc}
\VV_{1,1}^{\alpha} & \mtx{0} & \VV_{1,3}^{\alpha} \\
\mtx{0} & \VV_{2,2}^{\beta } & \VV_{2,3}^{\beta }  \\
\tfrac{1}{2}\VV_{3,1}^{\alpha} & \tfrac{1}{2}\VV_{3,2}^{\beta} & \tfrac{1}{2}\VV_{3,3}^{\alpha} + \tfrac{1}{2}\VV_{3,3}^{\beta}
\end{array}\right] +
\left[\begin{array}{c}
\VV_{1,4}^{\alpha} \\
\VV_{2,4}^{\beta}  \\
\tfrac{1}{2}\VV_{3,4}^{\alpha} + \tfrac{1}{2}\VV_{3,4}^{\beta}
\end{array}\right]\,\UU.
\end{equation}
An analogous computation for $\ww_{\rm e}$ yields
\begin{equation}
\label{eq:W_parent}
\WW =
\left[\begin{array}{ccc}
\WW_{1,1}^{\alpha} & \mtx{0} & \WW_{1,3}^{\alpha} \\
\mtx{0} & \WW_{2,2}^{\beta } & \WW_{2,3}^{\beta }  \\
\tfrac{1}{2}\WW_{3,1}^{\alpha} & \tfrac{1}{2}\WW_{3,2}^{\beta} & \tfrac{1}{2}\WW_{3,3}^{\alpha} + \tfrac{1}{2}\WW_{3,3}^{\beta}
\end{array}\right] +
\left[\begin{array}{c}
\WW_{1,4}^{\alpha} \\
\WW_{2,4}^{\beta}  \\
\tfrac{1}{2}\WW_{3,4}^{\alpha} + \tfrac{1}{2}\WW_{3,4}^{\beta}
\end{array}\right]\,\UU.
\end{equation}

\section{The full hierarchical scheme}
\label{sec:hierarchy}

\subsection{The algorithm}
\label{sec:thealgorithm}
Now that we know how to construct the DtN operator for a leaf (Section
\ref{sec:leaf}), and how to merge the DtN operators of two neighboring
patches to form the DtN operator of their union (Section
\ref{sec:merge}), we are ready to describe the full hierarchical
scheme for solving (\ref{eq:basic}).

\begin{figure}
\includegraphics[width=\textwidth]{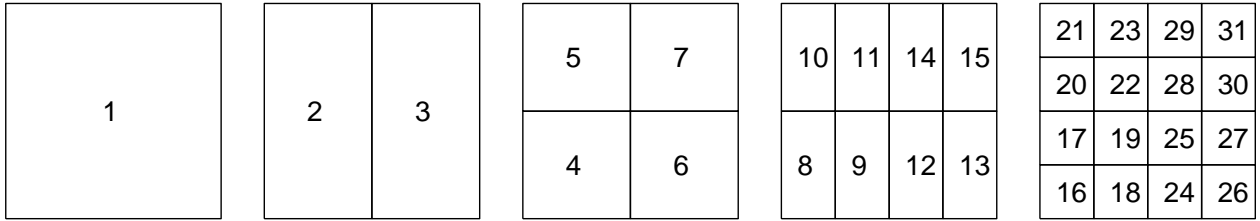}
\caption{The square domain $\Omega$ is split into $4 \times 4$ leaf boxes.
These are then gathered into a binary tree of successively larger boxes
as described in Section \ref{sec:thealgorithm}. One possible enumeration
of the boxes in the tree is shown, but note that the only restriction is
that if box $\tau$ is the parent of box $\sigma$, then $\tau < \sigma$.}
\label{fig:tree_numbering}
\end{figure}

First we partition the domain $\Omega$ into a collection of square (or
possibly rectangular) boxes, called \textit{leaf boxes}. These should
be small enough that a small spectral mesh with $p\times p$ nodes
(for, say, $p=20$) accurately interpolates both any potential solution
$u$ of (\ref{eq:basic}) and its partial derivatives $\partial_{1}u$,
$\partial_{2}u$, and $-\Delta u$. Let $\{\pxx_{k}\}_{k=1}^{N}$ denote
the points in this mesh. (Observe that nodes on internal boundaries
are shared between two or four local meshes.) Next construct a binary
tree on the collection of boxes by hierarchically merging them, making
sure that all boxes on the same level are roughly of the same size,
cf.~Figure \ref{fig:tree_numbering}.  The boxes should be ordered so
that if $\tau$ is a parent of a box $\sigma$, then $\tau < \sigma$. We
also assume that the root of the tree (i.e.~the full box $\Omega$) has
index $\tau=1$.

With each box $\tau$, we define two index vectors $I_{\rm i}^{\tau}$
and $I_{\rm e}^{\tau}$ as follows:

\vspace{1mm}

\begin{tabular}{ll}
$I_{\rm e}^{\tau}$ & A list of all indices of the nodes on the boundary of $\tau$.\\[2mm]

$I_{\rm i}^{\tau}$ & For a leaf $\tau$, $I_{\rm i}^{\tau}$ is a list of all interior nodes of $\tau$.\\
                & For a parent $\tau$, $I_{\rm i}^{\tau}$ is a list of all its interior nodes
                  that are not interior nodes of its children.
\end{tabular}

\vspace{1mm}

Next we execute a pre-computation in which for every box $\tau$, we construct the following matrices:

\vspace{1mm}

\begin{tabular}{ll}
$\UU^{\tau}$       & The matrix that maps the values of $\uu$ on the boundary of a box to the values of $\uu$ \\
                   & on the interior nodes of the box. In other words, $\uu(I_{\rm i}^{\tau}) = \UU^{\tau}\,\uu(I_{\rm e}^{\tau})$.\\[2mm]

$\VV^{\tau}$       & The matrix that maps the values of $\uu$ on the boundary of a box to the values of $\vv$\\
                   & (tabulating $du/dx_{1}$) on the boundary of a box.
                     In other words, $\vv(I_{\rm e}^{\tau}) = \VV^{\tau}\,\uu(I_{\rm e}^{\tau})$.\\[2mm]

$\WW^{\tau}$       & The matrix that maps the values of $\uu$ on the boundary of a box to the values of $\ww$\\
                   & (tabulating $du/dx_{2}$) on the boundary of a box.
                     In other words, $\ww(I_{\rm e}^{\tau}) = \WW^{\tau}\,\uu(I_{\rm e}^{\tau})$.
\end{tabular}

\vspace{1mm}

\noindent
To this end, we scan all boxes in the tree, going from smaller to larger.
For any leaf box $\tau$, a dense
matrix $\AA^{\tau}$ of size $p^{2}\times p^{2}$ that locally approximates the differential operator
in (\ref{eq:basic}) is formed. Then the matrices $\UU^{\tau}$, $\VV^{\tau}$, and $\WW^{\tau}$
are constructed via formulas (\ref{eq:U_leaf}), (\ref{eq:V_leaf}), and (\ref{eq:W_leaf}).
For a parent box $\tau$ with children $\sigma_{1}$ and $\sigma_{2}$, the matrices
$\UU^{\tau}$, $\VV^{\tau}$, and $\WW^{\tau}$ are formed from the DtN operators encoded in the
matrices $\VV^{\sigma_{1}}$, $\WW^{\sigma_{1}}$, $\VV^{\sigma_{2}}$, $\WW^{\sigma_{2}}$
using the formulas (\ref{eq:U_parent}), (\ref{eq:V_parent}), and (\ref{eq:W_parent}).
The full algorithm is summarized in Figure \ref{fig:precomp}. An illustrated cartoon
of the merge process is provided in Appendix \ref{app:cartoon}.

Once all the matrices $\{\UU^{\tau}\}_{\tau}$ have been formed, it is a simple matter
to construct a vector $\uu$ holding approximations to the solution $u$ of (\ref{eq:basic}).
The nodes are scanned starting with the root, and then proceeding down in the tree towards
smaller boxes. When a box $\tau$ is processed, the value of $\uu$ is known for all nodes
on its boundary (i.e.~those listed in $I_{\rm e}^{\tau}$). The matrix $\UU^{\tau}$ directly
maps these values to the values of $\uu$ on the nodes in the interior of $\tau$ (i.e.~those
listed in $I_{\rm i}^{\tau}$). When all nodes have been processed, all entries of $\uu$ have
been computed. Figure \ref{fig:solve} summarizes the solve stage.

\begin{remark}
Every interior meshpoint $\pxx_{k}$ belongs to the index vector $I_{\rm i}^{\tau}$
for precisely one node $\tau$. In other words $\bigcup_{\tau} I_{\rm i}^{\tau}$ forms a disjoint
union of the interior mesh points.
\end{remark}

\begin{remark}
\label{remark:tangential}
The way the algorithms are described, we compute for each node $\tau$ matrices
$\VV^{\tau}$ and $\WW^{\tau}$ that allow the computation of both the normal and
the tangential derivative at any boundary node, given the Dirichlet data on the
boundary. This is done for notational convenience only. In practice, any rows of
$\VV^{\tau}$ and $\WW^{\tau}$ that correspond to evaluation of tangential
derivatives need never be evaluated since tangential derivatives do not enter
into consideration at all.
\end{remark}


\begin{figure}
\fbox{
\begin{minipage}{140mm}
\begin{center}
\textsc{Pre-computation}
\end{center}

This program constructs the global Dirichlet-to-Neumann operator for (\ref{eq:basic}).\\
It also constructs all matrices $\UU^{\tau}$ required for constructing $u$ at all interior nodes.\\
It is assumed that if node $\tau$ is a parent of node $\sigma$, then $\tau < \sigma$.

\rule{\textwidth}{0.5pt}

\begin{tabbing}
\mbox{}\hspace{6mm} \= \mbox{}\hspace{6mm} \= \mbox{}\hspace{6mm} \= \mbox{}\hspace{6mm} \= \kill
\textbf{for} $\tau = N_{\rm boxes},\,N_{\rm boxes}-1,\,N_{\rm boxes}-2,\,\dots,\,1$\\
\> \textbf{if} ($\tau$ is a leaf)\\
\> \> $\vct{b}^{\tau}_{\rm i} = [b(\pxx_{j})]_{j \in I^{\tau}_{\rm i}}$\\
\> \> $\UU^{\tau} = -\bigl(-\mtx{L}_{\rm i,i} - \kappa^{2}\mbox{diag}(\vct{b}_{\rm i}^{\tau})\bigr)^{-1}\,\mtx{L}_{\rm i,e}$\\
\> \> $\VV^{\tau} = \mtx{D}_{\rm e,e} + \mtx{D}_{\rm e,i}\,\UU^{\tau}$\\
\> \> $\WW^{\tau} = \mtx{E}_{\rm e,e} + \mtx{D}_{\rm e,i}\,\UU^{\tau}$\\
\> \textbf{else}\\
\> \> Let $\sigma_{1}$ and $\sigma_{2}$ be the children of $\tau$.\\
\> \> Partition $I_{\rm e}^{\sigma_{1}}$ and $I_{\rm e}^{\sigma_{2}}$ into vectors $I_{1}$, $I_{2}$, $I_{3}$, and $I_{4}$ as shown in Figure \ref{fig:siblings_notation}.\\
\> \> \textbf{if} ($\sigma_{1}$ and $\sigma_{2}$ are side-by-side)\\
\> \> \> $\mtx{U}^{\tau} = \bigl( \VV^{\sigma_{1}}_{4,4} - \VV^{\sigma_{2}}_{4,4}\bigr)^{-1}
                           \bigl[-\VV^{\sigma_{1}}_{4,1}\  \big|\
                                  \VV^{\sigma_{2}}_{4,2}\  \big|\
                                  \VV^{\sigma_{2}}_{4,3} - \VV^{\sigma_{1}}_{4,3} \bigr]$\\
\> \> \textbf{else}\\
\> \> \> $\mtx{U}^{\tau} = \bigl( \WW^{\sigma_{1}}_{4,4} - \WW^{\sigma_{2}}_{4,4}\bigr)^{-1}
                           \bigl[-\WW^{\sigma_{1}}_{4,1}\  \big|\
                                  \WW^{\sigma_{2}}_{4,2}\  \big|\
                                  \WW^{\sigma_{2}}_{4,3} - \WW^{\sigma_{1}}_{4,3} \bigr]$\\
\> \> \textbf{end if}\\
\> \> $\VV^{\tau} = \left[\begin{array}{ccc}
                          \VV_{1,1}^{\sigma_{1}} & \mtx{0} & \VV_{1,3}^{\sigma_{1}} \\
                          \mtx{0} & \VV_{2,2}^{\sigma_{2} } & \VV_{2,3}^{\sigma_{2} }  \\
                          \tfrac{1}{2}\VV_{3,1}^{\sigma_{1}} & \tfrac{1}{2}\VV_{3,2}^{\sigma_{2}} &
                          \tfrac{1}{2}\VV_{3,3}^{\sigma_{1}} + \tfrac{1}{2}\VV_{3,3}^{\sigma_{2}}
                          \end{array}\right] +
                    \left[\begin{array}{c}
                          \VV_{1,4}^{\sigma_{1}} \\
                          \VV_{2,4}^{\sigma_{2}}  \\
                          \tfrac{1}{2}\VV_{3,4}^{\sigma_{1}} + \tfrac{1}{2}\VV_{3,4}^{\sigma_{2}}
                          \end{array}\right]\,\UU^{\tau}$.\\
\> \> $\WW^{\tau} = \left[\begin{array}{ccc}
                          \WW_{1,1}^{\sigma_{1}} & \mtx{0} & \WW_{1,3}^{\sigma_{1}} \\
                          \mtx{0} & \WW_{2,2}^{\sigma_{2} } & \WW_{2,3}^{\sigma_{2} }  \\
                          \tfrac{1}{2}\WW_{3,1}^{\sigma_{1}} & \tfrac{1}{2}\WW_{3,2}^{\sigma_{2}} &
                          \tfrac{1}{2}\WW_{3,3}^{\sigma_{1}} + \tfrac{1}{2}\WW_{3,3}^{\sigma_{2}}
                          \end{array}\right] +
                    \left[\begin{array}{c}
                          \WW_{1,4}^{\sigma_{1}} \\
                          \WW_{2,4}^{\sigma_{2}}  \\
                          \tfrac{1}{2}\WW_{3,4}^{\sigma_{1}} + \tfrac{1}{2}\WW_{3,4}^{\sigma_{2}}
                          \end{array}\right]\,\UU^{\tau}$.\\
\> \> Delete $\VV^{\sigma_{1}}$, $\WW^{\sigma_{1}}$, $\VV^{\sigma_{2}}$, $\WW^{\sigma_{2}}$.\\
\> \textbf{end if}\\
\textbf{end for}
\end{tabbing}
\end{minipage}}
\caption{Pre-computation}
\label{fig:precomp}
\end{figure}

\begin{figure}
\fbox{
\begin{minipage}{140mm}
\begin{center}
\textsc{Solver}
\end{center}

This program constructs an approximation $\uu$ to the solution $u$ of (\ref{eq:basic}).\\
It assumes that all matrices $\UU^{\tau}$ have already been constructed in a pre-computation.\\
It is assumed that if node $\tau$ is a parent of node $\sigma$, then $\tau < \sigma$.

\rule{\textwidth}{0.5pt}

\begin{tabbing}
\mbox{}\hspace{6mm} \= \mbox{}\hspace{6mm} \= \mbox{}\hspace{6mm} \= \mbox{}\hspace{6mm} \= \kill
$\uu(k) = f(\pxx_{k})$ for all $k \in I_{\rm e}^{1}$.\\
\textbf{for} $\tau = 1,\,2,\,3,\,\dots,\,N_{\rm boxes}$\\
\> $\uu(I_{\rm i}^{\tau}) = \UU^{\tau}\,\uu(I_{\rm i}^{\tau})$.\\
\textbf{end for}
\end{tabbing}
\end{minipage}}
\caption{Solve stage.}
\label{fig:solve}
\end{figure}

\begin{remark}
The merge stage is exact when performed in exact arithmetic. The only
approximation involved is the approximation of the solution $u$ on a leaf by its
interpolating polynomial.
\end{remark}

\subsection{Complexity analysis}
The analysis of the asymptotic cost of the algorithm in Section \ref{sec:thealgorithm}
closely mimics the analysis of the classical nested dissection algorithm
\cite{hoffman_1973,george_1973}. For simplicity, we analyze the simplest situation in
which a square domain is divided into $4^{L}$ leaf boxes, each holding a spectral
cartesian mesh with $p\times p$ points. The total number of unknowns in the system
is then roughly $4^{L}\,p^{2}$ (to be precise, $N = 4^{L}\,(p-1)^{2} + 2^{L+1}\,(p-1) + 1$).

\lsp

\noindent
\textit{Cost of leaf computation:} Evaluating the formulas
(\ref{eq:U_leaf}), (\ref{eq:V_leaf}), and (\ref{eq:W_leaf}) requires dense
matrix algebra on matrices of size roughly $p^{2} \times p^{2}$. Since there
are about $N/p^{2}$ leaves, the total cost is
$$
T_{\rm leaf} \sim \frac{N}{p^{2}} \times \left(p^{2}\right)^{3} \sim N\,p^{4}.
$$

\lsp

\noindent
\textit{Cost of the merge operations:} For an integer $\ell \in \{0,\,1,\,2,\,\dots,\,L\}$,
we refer to ``level $\ell$'' as the collection of boxes whose side length is $2^{-\ell}$ times
the side of the full box $\Omega$ (so that $\ell=0$ corresponds to the root and $\ell=L$ corresponds
to the set of leaf boxes). To form the matrices $\UU^{\tau}$, $\VV^{\tau}$, and $\WW^{\tau}$ for a
box on level $\ell$, we need to evaluate each of the formulas (\ref{eq:U_parent}), (\ref{eq:V_parent}),
and (\ref{eq:W_parent}) three times, with each computation involving matrices of size roughly
$2^{-\ell}N^{0.5} \times 2^{-\ell}N^{0.5}$. Since there are $4^{\ell}$ boxes on level $\ell$,
we find that the cost of processing level $\ell$ is
$$
T_{\ell} \sim 4^{\ell} \times \left(2^{-\ell}\,N^{0.5}\right)^{3} \sim 2^{-\ell}\,N^{1.5}.
$$
Adding the costs at all levels, we get
$$
T_{\rm merge} \sim
\sum_{\ell=0}^{L-1} T_{\ell} \sim
\sum_{\ell=0}^{L-1} 2^{-\ell}\,N^{1.5} \sim
N^{1.5}.
$$

\lsp

\noindent
\textit{Cost of solve stage:} The cost of processing a non-leaf node on level $\ell$ is simply the
cost of a matrix-vector multiply involving the dense matrix $\UU^{\tau}$ of size $2^{-\ell}N^{0.5} \times 2^{-\ell}N^{0.5}$.
For a leaf, $\UU^{\tau}$ is of size roughly $p^{2} \times p$. Therefore
$$
T_{\rm leaf} \sim \frac{N}{p^{2}} \times p^{2}\,p + \sum_{\ell=0}^{L-1}4^{\ell} \times \left(2^{-\ell}N^{0.5}\right)^{2} \sim
N\,p + N\,L \sim N\,p + N\,\log(N).
$$

\subsection{Problem of resonances}
\label{sec:resonance}

The scheme presented in Section \ref{sec:thealgorithm} will fail if one of the patches
in the hierarchical partitioning is resonant in the sense that there exist non-trivial
solutions to the Helmholtz equation that have zero Dirichlet data at the boundary of
the patch. In this case, the Neumann data for the patch is not uniquely determined by the
Dirichlet data, and the DtN operator cannot exist. In practice, this problem will of course
arise even if we are merely \textit{close} to a resonance, and will be detected by
the discovery that the inverse matrix in the formulas (\ref{eq:U_leaf}) and (\ref{eq:U_parent})
for the solution operator $\UU^{\tau}$ is ill-conditioned.

It is our experience from working with domains of size one hundred wave-lengths or less
that resonances are very rare; one almost never encounters the problem.
Nevertheless, it is important to monitor the conditioning of the formulas (\ref{eq:U_leaf})
and (\ref{eq:U_parent}) to ensure the accuracy of the final answer. Should a problem be
detected, the easiest solution would be to simply start the computation over with a
different tessellation of the domain $\Omega$. Very likely, this will resolve the problem.

Current efforts to formulate variations of the scheme that are inherently not
vulnerable to resonances are described in Section \ref{sec:reformulation}.


\section{Numerical experiments}
\label{sec:num}

This section reports the results of some numerical experiments with
the method described in Section \ref{sec:thealgorithm}. The method was
implemented in Matlab and the experiments executed on a Lenovo W510
laptop with a quad core Intel i7 Q720 processor with 1.6GHz
clockspeed, and 16GB of RAM.

The speed and memory requirements of the algorithm were investigated
by solving the special case where $b=0$ in (\ref{eq:basic}), and the
Dirichlet data is simply set to equal a known analytic solution.
The results are reported in Section
\ref{sec:helmholtz_special}. We also report the errors incurred in the
special case, but it should be noted that these represent only a best
case estimate of the errors since the equation solved is particularly benign.
To get a more realistic estimate of the errors in the
method, we also applied it to three situation in which exact solutions
are not known: A problem with variable coefficients in Section
\ref{sec:helmholtz_bump}, a problem on an L-shaped domain in Section
\ref{sec:L}, and a convection-diffusion problem in Section
\ref{sec:convdiff}.

\subsection{Constant coefficient Helmholtz}
\label{sec:helmholtz_special}
We solved the basic Helmholtz equation
\begin{align}
-\Delta u(\pxx) - \kappa^{2}\,u(\pxx) =&\ 0       \qquad& \pxx \in \Omega,\\
                              u(\pxx) =&\ f(\pxx) \qquad& \pxx \in \Gamma,
\end{align}
where $\Omega = [0,1]^{2}$ and $\Gamma = \partial \Omega$. The boundary data were
in this first set of experiments chosen as the restriction to $\Gamma$ of an exact
solution
\begin{equation}
\label{eq:uexact}
u_{\rm exact}(\pxx) = Y_{0}(\kappa|\pxx - \hat{\pxx}|),
\end{equation}
where $\hat{\pxx} = (-0.2,\,0.4)$, and where $Y_{0}$ is the $0$'th
Bessel function of the second kind.  This experiment serves two
purposes. The first is to systematically measure the speed and memory
requirements of the method described in Section
\ref{sec:thealgorithm}. The second is to get a sense of what errors
can be expected in a ``best case'' scenario with a very smooth
solution.  Observe however that the situation is by no means
artificial since the smoothness of this case is exactly what one
encounters when the solver is applied to a free space scattering
problem as described in Section \ref{sec:freespace}.

The domain $\Omega$ was discretized into $n\times n$ patches, and on
each patch a $p\times p$ Cartesian mesh of Chebyshev nodes was
placed. The total number of degrees of freedom is then
$$
N = \bigl(n(p-1)\bigr)^{2} + 2n(p-1) + 1.
$$ We tested the method for $p \in \{6,\,11,\,21,\,41\}$. For each
fixed $p$, the method was executed for several different mesh sizes
$n$. The wave-number $\kappa$ was chosen to keep a constant of 12
points per wavelength, or $\kappa = 2\pi\,n(p-1)/12$.

Since the exact solution is known in this case, we computed the direct error measure
$$
E_{\rm pot} = \max_{k\,\colon\,\pxx_{k} \in \Omega} \bigl|u_{\rm computed}(\pxx_{k}) - u_{\rm exact}(\pxx_{k})\bigr|,
$$
where $\{\pxx_{k}\}_{k=1}^{N}$ is the set of all mesh points. We also computed the
maximum error in the gradient of $u$ on the boundary as computed via the $\VV$ and $\WW$
operators on the root box,
$$
E_{\rm grad} = \max_{k\,\colon\,\pxx_{k} \in \Gamma}
\bigl\{\bigl|v_{\rm computed}(\pxx_{k}) - [\partial_{1}u_{\rm exact}](\pxx_{k})\bigr|,\,
       \bigl|w_{\rm computed}(\pxx_{k}) - [\partial_{2}u_{\rm exact}](\pxx_{k})\bigr|\bigr\}.
$$

Table \ref{table:helmholtz_special_12} reports the following variables:\\
\begin{tabular}{ll}
$N_{\rm wave}$  & The number of wave-lengths along one side of $\Omega$.\\
$t_{\rm inv}$   & The time in seconds required to execute the pre-computation in Figure \ref{fig:precomp}.\\
$t_{\rm solve}$ & The time in seconds required to execute the solve in Figure \ref{fig:solve}.\\
$M$             & The amount of RAM used in the pre-computation in MB.
\end{tabular}\\
The table also reports the memory requirements in terms of number of the number of double precision reals
that need to be stored per degree of freedom in the discretization.

The high-order version of the method ($p=41$) was also capable of performing a high accuracy solve
with only six points per wave-length. The results are reported in Table \ref{table:helmholtz_special_6}.

We see that increasing the spectral order is very beneficial for improving accuracy. However, the
speed deteriorates and the memory requirements increase as $p$ grows. Choosing $p=21$ appears to
be a good compromise.

%

\begin{remark}
In the course of executing the numerical examples, the instability
problem described in Section \ref{sec:resonance} was detected
precisely once (for $p=16$ and 12 points per wave length).
\end{remark}

\begin{table}
\small
\begin{tabular}{|r|r|r|r|r|r|r|r|r|r|r|r|}\hline
$p$ & $N$ & $N_{\rm wave}$ & $t_{\rm inv}$ & $t_{\rm solve}$ & $E_{\rm pot}$ & $E_{\rm grad}$ & $M$ & $M/N$ \\
&&& (sec) & (sec) &&& (MB) & (reals/DOF) \\ \hline
 6 &    6561 &   6.7 &    0.28 &    0.0047 &  8.02105e-03 &  3.06565e-01 &       2.8 &         56.5 \\
 6 &   25921 &  13.3 &    0.96 &    0.0184 &  1.67443e-02 &  1.33562e+00 &      12.7 &         64.2 \\
 6 &  103041 &  26.7 &    4.42 &    0.0677 &  3.60825e-02 &  5.46387e+00 &      56.2 &         71.5 \\
 6 &  410881 &  53.3 &   20.23 &    0.2397 &  3.39011e-02 &  1.05000e+01 &     246.9 &         78.8 \\
 6 & 1640961 & 106.7 &   88.73 &    0.9267 &  7.48385e-01 &  4.92943e+02 &    1075.0 &         85.9 \\ \hline
11 &    6561 &   6.7 &    0.16 &    0.0019 &  2.67089e-05 &  1.08301e-03 &       2.9 &         58.0 \\
11 &   25921 &  13.3 &    0.68 &    0.0073 &  5.30924e-05 &  4.34070e-03 &      13.0 &         65.7 \\
11 &  103041 &  26.7 &    3.07 &    0.0293 &  1.01934e-04 &  1.60067e-02 &      57.4 &         73.0 \\
11 &  410881 &  53.3 &   14.68 &    0.1107 &  1.07747e-04 &  3.49637e-02 &     251.6 &         80.2 \\
11 & 1640961 & 106.7 &   68.02 &    0.3714 &  2.17614e-04 &  1.37638e-01 &    1093.7 &         87.4 \\ \hline
21 &    6561 &   6.7 &    0.23 &    0.0011 &  2.56528e-10 &  1.01490e-08 &       4.4 &         87.1 \\
21 &   25921 &  13.3 &    0.92 &    0.0044 &  5.24706e-10 &  4.44184e-08 &      18.8 &         95.2 \\
21 &  103041 &  26.7 &    4.68 &    0.0173 &  9.49460e-10 &  1.56699e-07 &      80.8 &        102.7 \\
21 &  410881 &  53.3 &   22.29 &    0.0727 &  1.21769e-09 &  3.99051e-07 &     344.9 &        110.0 \\
21 & 1640961 & 106.7 &   99.20 &    0.2965 &  1.90502e-09 &  1.24859e-06 &    1467.2 &        117.2 \\
21 & 6558721 & 213.3 &  551.32 &   20.9551 &  2.84554e-09 &  3.74616e-06 &    6218.7 &        124.3 \\ \hline
41 &    6561 &   6.7 &    1.50 &    0.0025 &  9.88931e-14 &  3.46762e-12 &       7.9 &        157.5 \\
41 &   25921 &  13.3 &    4.81 &    0.0041 &  1.58873e-13 &  1.12883e-11 &      32.9 &        166.4 \\
41 &  103041 &  26.7 &   18.34 &    0.0162 &  3.95531e-13 &  5.51141e-11 &     137.1 &        174.4 \\
41 &  410881 &  53.3 &   75.78 &    0.0672 &  3.89079e-13 &  1.03546e-10 &     570.2 &        181.9 \\
41 & 1640961 & 106.7 &  332.12 &    0.2796 &  1.27317e-12 &  7.08201e-10 &    2368.3 &        189.2 \\ \hline
\end{tabular}
\caption{Results from an experiment with a constant coefficient Helmholtz problem on a square. The
boundary data were picked so that the analytic solution was known; as a consequence, the solution
is smooth, and can be smoothly extended across the boundary. The wave-number was chosen to keep a
constant of 12 discretization points per wave-length.}
\label{table:helmholtz_special_12}
\end{table}

\begin{table}
\small
\begin{tabular}{|r|r|r|r|r|r|r|r|r|r|r|r|}\hline
$p$ & $N$ & $N_{\rm wave}$ & $t_{\rm inv}$ & $t_{\rm solve}$ & $E_{\rm pot}$ & $E_{\rm grad}$ & $M$ & $M/N$\\ \hline
&&& (sec) & (sec) &&& (MB) & (reals/DOF) \\ \hline
41 &    6561 &  13.3 &    1.30 &    0.0027 &  1.54407e-09 &  1.78814e-07 &       7.9 &        157.5 \\
41 &   25921 &  26.7 &    4.40 &    0.0043 &  1.42312e-08 &  2.35695e-06 &      32.9 &        166.4 \\
41 &  103041 &  53.3 &   17.54 &    0.0199 &  1.73682e-08 &  5.84193e-06 &     137.1 &        174.4 \\
41 &  410881 & 106.7 &   72.90 &    0.0717 &  2.28475e-08 &  1.51575e-05 &     570.2 &        181.9 \\
41 & 1640961 & 213.3 &  307.37 &    0.3033 &  4.12809e-08 &  5.51276e-05 &    2368.3 &        189.2 \\ \hline
\end{tabular}
\caption{This table illustrates the same situation as Table \ref{table:helmholtz_special_12},
but now $\kappa$ is increased twice as fast (so that we keep only 6 points per wave-length).}
\label{table:helmholtz_special_6}
\end{table}

\subsection{Variable coefficient Helmholtz}
\label{sec:helmholtz_bump}
We solved the equation
\begin{align}
\label{eq:helmholtz_bump}
-\Delta u(\pxx) - \kappa^{2}\,\bigl(1 - b(\pxx)\bigr)\,u(\pxx) =&\ 0       \qquad& \pxx \in \Omega,\\
                                       u(\pxx) =&\ f(\pxx) \qquad& \pxx \in \Gamma,
\end{align}
where $\Omega = [0,1]^{2}$, where $\Gamma = \partial \Omega$, and where
$$
b(\pxx) = \left(\sin(4\pi x_{1})\,\sin(4\pi x_{2})\right)^{2}.
$$
The Helmholtz parameter was kept fixed at $\kappa=80$, corresponding to a domain
size of $12.7 \times 12.7$ wave lengths. The boundary data was given by
$$
f(\pxx) = \cos(8 x_{1})\,\bigl(1 - 2x_{2}\bigr).
$$
The equation (\ref{eq:helmholtz_bump}) was discretized and solved as described
in Section \ref{sec:helmholtz_special}. The speed and memory requirements for
this computation are exactly the same as for the example in Section \ref{sec:helmholtz_special}
(they do not depend on what equation is being solved), so we now focus on the accuracy
of the method. We do not know of an exact solution, and therefore report
pointwise convergence. Letting $u_{N}$ denote the value of $u$ computed using $N$ degrees
of freedom, we used
$$
E_{N}^{\rm int} = u_{N}(\hat{\pxx}) - u_{4N}(\hat{\pxx})
$$
as an estimate for the pointwise error in $u$ at the point $\hat{\pxx} = (0.75,\,0.25)$.
We analogously estimated convergence of the normal derivative at the point
$\hat{\pyy} = (0.75,\,0.00)$ by measuring
$$
E_{N}^{\rm bnd} = w_{N}(\hat{\pyy}) - w_{4N}(\hat{\pyy}).
$$
The results are reported in Table \ref{table:helmholtz_bump12}.
Table \ref{table:helmholtz_bump102} reports the results from an
analogous experiment, but now for a domain of size $102 \times 102$ wave-lengths.

We observe that accuracy is almost as good as for the constant
coefficient case.  Ten digits of accuracy is easily attained, but
getting more than that seems challenging; increasing $N$ further leads
to no improvement in accuracy. The method appears to be stable in the
sense that nothing bad happens when $N$ is either too large or too
small.

\begin{table}
\small
\begin{tabular}{|r|r|r|r|r|r|r|r|r|r|r|r|r|r|r|} \hline
$p$ & $N$ & pts per wave & $u_{N}(\hat{\pxx})$ & $E_{N}^{\rm int}$ & $w_{N}(\hat{\pyy})$ & $E_{N}^{\rm bnd}$ \\ \hline
 6 &    6561 &    6.28 &  -2.505791196753718 & -2.457e-01 &    -661.0588680825 & -8.588e+03 \\
 6 &   25921 &   12.57 &  -2.260084219562163 &  1.676e-01 &    7926.8096554095 &  8.141e+04 \\
 6 &  103041 &   25.13 &  -2.427668162910011 &  1.779e-02 &  -73484.9989261573 & -3.894e+04 \\
 6 &  410881 &   50.27 &  -2.445455646843485 &  1.233e-03 &  -34547.6403539568 & -1.235e+03 \\
 6 & 1640961 &  100.53 &  -2.446688310709834 &  7.891e-05 &  -33313.0000081604 & -7.627e+01 \\
 6 & 6558721 &  201.06 &  -2.446767218259172 &            &  -33236.7252190062 &            \\ \hline
11 &    6561 &    6.28 &  -2.500353149793093 & -5.375e-02 &  -27023.0713474340 &  6.524e+03 \\
11 &   25921 &   12.57 &  -2.446599788642489 &  1.728e-04 &  -33547.3621639994 & -3.153e+02 \\
11 &  103041 &   25.13 &  -2.446772604281610 & -9.465e-08 &  -33232.0940315585 & -4.754e-01 \\
11 &  410881 &   50.27 &  -2.446772509631734 &  3.631e-10 &  -33231.6186528531 & -7.331e-04 \\
11 & 1640961 &  100.53 &  -2.446772509994819 &            &  -33231.6179197169 &            \\ \hline
21 &    6561 &    6.28 &  -2.448236804078803 & -1.464e-03 &  -32991.4583727724 &  2.402e+02 \\
21 &   25921 &   12.57 &  -2.446772430608166 &  7.976e-08 &  -33231.6118304666 &  5.984e-03 \\
21 &  103041 &   25.13 &  -2.446772510369452 &  5.893e-11 &  -33231.6178142514 & -5.463e-06 \\
21 &  410881 &   50.27 &  -2.446772510428384 &  2.957e-10 &  -33231.6178087887 & -2.792e-05 \\
21 & 1640961 &  100.53 &  -2.446772510724068 &            &  -33231.6177808723 &            \\ \hline
41 &    6561 &    6.28 &  -2.446803898373796 & -3.139e-05 &  -33233.0037457220 & -1.386e+00 \\
41 &   25921 &   12.57 &  -2.446772510320572 &  1.234e-10 &  -33231.6179029824 & -8.940e-05 \\
41 &  103041 &   25.13 &  -2.446772510443995 &  2.888e-11 &  -33231.6178135860 & -1.273e-05 \\
41 &  410881 &   50.27 &  -2.446772510472872 &  7.731e-11 &  -33231.6178008533 & -4.668e-05 \\
41 & 1640961 &  100.53 &  -2.446772510550181 &            &  -33231.6177541722 &            \\ \hline
\end{tabular}
\caption{Results from a variable coefficient Helmholtz problem on a domain of size $12.7 \times 12.7$ wave-lengths.}
\label{table:helmholtz_bump12}
\end{table}

\begin{table}
\small
\begin{tabular}{|r|r|r|r|r|r|r|r|r|r|r|r|r|r|r|} \hline
$p$ & $N$ & pts per wave & $u_{N}(\hat{\pxx})$ & $E_{N}^{\rm int}$ & $w_{N}(\hat{\pyy})$ & $E_{N}^{\rm bnd}$ \\ \hline
21 &    6561 &    0.79 &   0.007680026148649 &  4.085e-03 &   3828.84075823538 &  6.659e+03 \\
21 &   25921 &    1.57 &   0.003595286353011 &  1.615e+00 &  -2829.88055527014 & -1.791e+02 \\
21 &  103041 &    3.14 &  -1.611350573683137 &  1.452e+00 &  -2650.80640712917 & -5.951e+03 \\
21 &  410881 &    6.28 &  -3.063762877533994 &  4.557e-03 &   3299.72573600854 & -7.772e+00 \\
21 & 1640961 &   12.57 &  -3.068320356836451 & -7.074e-08 &   3307.49786015114 &  9.592e-05 \\
21 & 6558721 &   25.13 &  -3.068320286093162 &            &   3307.49776422768 &            \\ \hline
41 &    6561 &    0.79 &  -0.000213617359480 & -1.608e-01 &  -833.919575889393 & -1.006e+03 \\
41 &   25921 &    1.57 &   0.160581838547352 & -7.415e-01 &   171.937456004515 & -1.797e+03 \\
41 &  103041 &    3.14 &   0.902057033817060 &  3.970e+00 &  1969.187023322940 & -1.338e+03 \\
41 &  410881 &    6.28 &  -3.068320045777766 &  2.405e-07 &  3307.497852217234 &  8.792e-05 \\
41 & 1640961 &   12.57 &  -3.068320286282830 &(-1.897e-10)&  3307.497764294187 & (6.651e-8) \\ \hline
\end{tabular}
\caption{Results from a variable coefficient Helmholtz problem on a domain of size $102 \times 102$ wave-lengths.}
\label{table:helmholtz_bump102}
\end{table}

\subsection{L-shaped domain}
\label{sec:L}
We solved the equation
\begin{align}
-\Delta u(\pxx) - \kappa^{2}\,u(\pxx) =&\ 0       \qquad& \pxx \in \Omega,\\
                              u(\pxx) =&\ f(\pxx) \qquad& \pxx \in \Gamma,
\end{align}
where $\Omega$ is the L-shaped domain
$$
\Omega = [0,2]^{2} \backslash [1,2]^{2},
$$
and where the Helmholtz parameter $\kappa$ is held fixed at $\kappa = 40$,
making the domain $12.7 \times 12.7$ wave-lengths large. The pointwise
errors were estimated at the points
$$
\hat{\pxx} = (0.75,\,0.75),
\qquad\mbox{and}\qquad
\hat{\pyy} = (1.25,\,1.00),
$$
via
$$
E_{N}^{\rm int} = u_{N}(\hat{\pxx}) - u_{4N}(\hat{\pxx}),
\qquad\mbox{and}\qquad
E_{N}^{\rm bnd} = w_{N}(\hat{\pxx}) - w_{4N}(\hat{\pxx}).
$$
The results are given in Table \ref{table:helmholtz_L}.

We observe that the errors are in this case significantly larger than they were
for square domains. This is presumably due to the fact that the solution $u$ is
singular near the re-entrant corner in the L-shaped domain. (When boundary conditions
corresponding to an exact solution like (\ref{eq:uexact}) are imposed, the method is
just as accurate as it is for the square domain.) Nevertheless, the method easily
attains solutions with between four and five correct digits.

\begin{table}
\small
\begin{tabular}{|r|r|r|r|r|r|r|r|r|r|r|r|r|r|r|} \hline
$p$ & $N$ & pts per wave & $u_{N}(\hat{\pxx})$ & $E_{N}^{\rm int}$ & $w_{N}(\hat{\pyy})$ & $E_{N}^{\rm bnd}$ \\ \hline
 6 &   19602 & 12.57 &  8.969213152495405 &  2.258e+00 &   226.603823940515 &  8.748e+01 \\
 6 &   77602 & 25.13 &  6.711091204119065 &  1.317e-01 &   139.118986915759 &  4.949e+00 \\
 6 &  308802 & 50.27 &  6.579341284597024 &  8.652e-03 &   134.169908546083 &  3.261e-01 \\
 6 & 1232002 &100.53 &  6.570688999911585 &            &   133.843774958376 &            \\ \hline
11 &   19602 & 12.57 &  6.571117172871830 &  9.613e-04 &   133.865552472382 &  3.851e-02 \\
11 &   77602 & 25.13 &  6.570155895761215 &  5.154e-05 &   133.827043929015 &  5.207e-03 \\
11 &  308802 & 50.27 &  6.570104356719250 &  1.987e-05 &   133.821836691967 &  2.052e-03 \\
11 & 1232002 &100.53 &  6.570084491282650 &            &   133.819785089497 &            \\ \hline
21 &   19602 & 12.57 &  6.570152809642857 &  4.905e-05 &   133.898328735897 &  7.663e-02 \\
21 &   77602 & 25.13 &  6.570103763348836 &  1.951e-05 &   133.821703687416 &  1.943e-03 \\
21 &  308802 & 50.27 &  6.570084254517955 &  7.743e-06 &   133.819760759394 &  7.996e-04 \\
21 & 1232002 &100.53 &  6.570076511737839 &            &   133.818961147570 &            \\ \hline
\end{tabular}
\caption{Results from a constant coefficient Helmholtz problem on an L-shaped domain of size $12.7 \times 12.7$ wave-lengths.}
\label{table:helmholtz_L}
\end{table}

\subsection{Convection diffusion}
\label{sec:convdiff}
We solved the equation
\begin{align}
\label{eq:convdiff}
-\Delta u(\pxx) - 1000\,[\partial_{2}u](\pxx)  =&\ 0       \qquad& \pxx \in \Omega,\\
                                       u(\pxx) =&\ f(\pxx) \qquad& \pxx \in \Gamma,
\end{align}
where $\Omega = [0,1]^{2}$, where $\Gamma = \partial \Omega$, and where
the boundary data was given by
$$
f(\pxx) = \cos(x_{1})\,e^{x_{2}}.
$$
The equation (\ref{eq:helmholtz_bump}) was discretized and solved as described
in Section \ref{sec:helmholtz_special}. Note that this case involves a non-oscillatory
solution and does not fit the template (\ref{eq:basic}). It is included to illustrate
how the spectral method handles a sharp gradient in the solution.

The pointwise errors were estimated at the points
$$
\hat{\pxx} = (0.75,\,0.25),
\qquad\mbox{and}\qquad
\hat{\pyy} = (0.75,\,0.00),
$$
via
$$
E_{N}^{\rm int} = u_{N}(\hat{\pxx}) - u_{4N}(\hat{\pxx}),
\qquad\mbox{and}\qquad
E_{N}^{\rm bnd} = w_{N}(\hat{\pxx}) - w_{4N}(\hat{\pxx}).
$$
The results are given in Table \ref{table:convdiff1000}.
Table \ref{table:convdiff10000} reports results from an analogous experiment,
but now with the strength of the convection term further increased by a factor of 10.

We observe that the method has no difficulties resolving steep gradients,
and that moderate order methods ($p=11$) perform very well here.

\begin{table}
\small
\begin{tabular}{|r|r|r|r|r|r|r|r|r|r|r|r|r|r|r|} \hline
$p$ & $N$ & $u_{N}(\hat{\pxx})$ & $E_{N}^{\rm int}$ & $w_{N}(\hat{\pyy})$ & $E_{N}^{\rm bnd}$ \\ \hline
11 &   25921 &   1.987126286905920 & -3.191e-04 &  1255.25512379751 & -7.191e-03 \\
11 &  103041 &   1.987445414657945 &  3.979e-13 &  1255.26231503666 & -6.529e-04 \\
11 &  410881 &   1.987445414657547 &  2.455e-12 &  1255.26296795281 & -1.889e-05 \\
11 & 1640961 &   1.987445414655092 &            &  1255.26298684450 &            \\ \hline
21 &   25921 &   1.987076984861468 & -3.684e-04 &  1255.26075989546 & -2.186e-03 \\
21 &  103041 &   1.987445414658047 & -3.009e-13 &  1255.26294637880 & -4.054e-05 \\
21 &  410881 &   1.987445414658348 & -2.600e-13 &  1255.26298691798 & -7.881e-08 \\
21 & 1640961 &   1.987445414658608 &            &  1255.26298699680 &            \\ \hline
41 &   25921 &   1.988004762686629 &  5.593e-04 &  1255.26290210213 & -8.478e-05 \\
41 &  103041 &   1.987445414657579 & -9.706e-13 &  1255.26298687891 & -1.178e-07 \\
41 &  410881 &   1.987445414658550 & -1.237e-12 &  1255.26298699669 & -1.636e-09 \\
41 & 1640961 &   1.987445414659787 &            &  1255.26298699832 &            \\ \hline
\end{tabular}
\caption{Errors for the convection diffusion problem (\ref{eq:convdiff}).}
\label{table:convdiff1000}.
\end{table}

\begin{table}
\small
\begin{tabular}{|r|r|r|r|r|r|r|r|r|r|r|r|r|r|r|} \hline
$p$ & $N$ & $u_{N}(\hat{\pxx})$ & $E_{N}^{\rm int}$ & $w_{N}(\hat{\pyy})$ & $E_{N}^{\rm bnd}$ \\ \hline
21 &   25921 &   1.476688750775769 & -4.700e-01 &  13002.9937044202 &  4.325e+02 \\
21 &  103041 &   1.946729131937971 & -4.206e-02 &  12570.4750256324 & -7.862e-03 \\
21 &  410881 &   1.988785675941193 & -1.716e-06 &  12570.4828877374 & -4.900e-03 \\
21 & 1640961 &   1.988787391699051 & (6.719e-13)&  12570.4877875310 &(-4.411e-04)\\ \hline
41 &   25921 &   2.587008191566030 &  6.407e-01 &  13002.1084152522 &  4.316e+02 \\
41 &  103041 &   1.946284950165041 & -4.250e-02 &  12570.4835546978 & -2.618e-03 \\
41 &  410881 &   1.988785277235741 & -2.114e-06 &  12570.4861729647 & -2.127e-03 \\
41 & 1640961 &   1.988787391699218 &            &  12570.4882994934 &            \\ \hline
\end{tabular}
\caption{Errors for a convection diffusion problem similar to (\ref{eq:convdiff}), but now
for the even more convection dominated operator $A = -\Delta - 10\,000\,\partial_{2}$.}
\label{table:convdiff10000}
\end{table}


\section{Extensions}
\label{sec:extensions}

\subsection{Linear complexity algorithms}
\label{sec:linear}
In discussing the asymptotic complexity of the scheme in Section \ref{sec:thealgorithm}
it is important to distinguish between the case where $N$ is increased for a fixed
wave-number $\kappa$, and the case where $\kappa \sim N^{0.5}$ to keep the number of
degrees of freedom per wavelength constant.

Let us first discuss the case where the wave-number $\kappa$ is kept fixed as $N$ is
increased. In this situation, the matrices $\UU^{\tau}$ will for the parent nodes
become highly rank deficient (to finite precision). By factoring these matrices in
the pre-computation, the solve phase can be reduced to $O(N)$ complexity. Moreover,
the off-diagonal blocks of the matrices $\VV^{\tau}$ and $\WW^{\tau}$ will also be
of low numerical rank; technically, they can efficiently be represented in data sparse
formats such as the $\mathcal{H}$-matrix format \cite{2003_hackbusch}, or the
\textit{Hierarchically Semi-Separable}-format of \cite{2012_martinsson_FDS_survey,2010_gu_xia_HSS}.
This property can be exploited to reduce the complexity of the
pre-computation from $O(N^{1.5})$ to $O(N)$ in a manner similar to what is done for
classical nested dissection for finite-element and finite-difference matrices in
\cite{2011_gillman_dissertation,2007_leborne_HLU,2011_ying_nested_dissection_2D,2009_xia_superfast}.
Note that while the acceleration of the solve phase is trivial, it takes some work
to exploit the more complicated structure in $\VV^{\tau}$ and $\WW^{\tau}$.

The case where $\kappa\sim N^{0.5}$ as $N$ increases is more complicated. In this situation,
the numerical ranks of the matrices $\UU^{\tau}$ and the off-diagonal blocks of $\VV^{\tau}$
and $\WW^{\tau}$ will be large enough that no reduction in asymptotic complexity can be
expected. However, the matrices are in practice far from full rank, and substantial savings
can be achieved by exploiting techniques such as those described in the previous paragraph.

\subsection{General elliptic problems}
The algorithm in Section \ref{sec:thealgorithm} can straight-forwardly
be generalized to elliptic equations like
\begin{multline*}
-c_{11}(\pxx)[\partial_{1}^{2}u](\pxx)
-2c_{12}(\pxx)[\partial_{1}\partial_{2}u](\pxx)
-c_{22}(\pxx)[\partial_{2}^{2}u](\pxx)\\
+c_{1}(\pxx)[\partial_{1}u](\pxx)
+c_{2}(\pxx)[\partial_{2}u](\pxx)
+c(\pxx)\,u(\pxx) = 0
\end{multline*}
as long as the coefficients are smooth and the leading order operator is elliptic.
The only modification required is that in the leaf computation, the matrix
representing a spectral approximation of the differential operator be replaced by
something like
$$
\mtx{A}^{\tau} =
-\mtx{C}_{11}\mtx{D}^{2}
-2\mtx{C}_{12}\mtx{D}\mtx{E}
-\mtx{C}_{22}\mtx{E}^{2}
+\mtx{C}_{1}\mtx{D}
+\mtx{C}_{2}\mtx{E}
+\mtx{C},
$$
where $\mtx{C}_{11}$ is a diagonal matrix with entries $c_{11}(\pvct{x}_{j})$ for $j\in I^{\tau}$, etc.
Observe that only the leaf computation needs to be modified, the merge steps remain exactly the same.
The stability and accuracy of the method of course depends on the signs and
relative magnitudes of the coefficient terms, but tentative numerical experiments
indicate that the method is effective for a broad range of different problems,
including convection-dominated convection-diffusion equations.

\subsection{Free space scattering problems in the plane}
\label{sec:freespace}
If you know the DtN operator for the inhomogeneous square, you can very rapidly
solve an exterior scattering problem as follows. Consider the equation
$$
-\Delta u(\pxx) - \kappa^{2}\,(1 - b(\pxx))\,u(\pxx) = f(\pxx),\qquad \pxx \in \mathbb{R}^{2}.
$$
Appropriate radiation conditions at infinity are imposed on $u$.
We assume that $b$ is compactly supported inside the domain $\Omega$, and that
$f$ is supported outside $\Omega$. The standard technique is to look for a solution
$u$ of the form
$$
u = v + w,
$$
where $v$ is an \textit{incoming field} and $w$ is the \textit{outgoing field}.
The incoming field is defined by
\begin{equation}
\label{eq:v}
v(\pxx) = [\phi_{\kappa} * f](\pxx) = \int_{\mathbb{R}^{2}}\phi_{\kappa}(\pxx-\pyy)\,f(\pyy)\,d\pyy,
\end{equation}
where $\phi_{\kappa}$ is the fundamental solution to the free space Helmholtz problem.
Then $-\Delta v - \kappa^{2}v = f$, and since $b(\pxx) = 0$ for $\pxx \in \Omega^{\rm c}$,
we find that the outgoing potential $w$ must satisfy
\begin{equation}
\label{eq:comstock2}
-\Delta w(\pxx) - \kappa^{2}w(\pxx) = 0,\qquad \pxx \in \Omega^{\rm c}.
\end{equation}
Now use the method in Section \ref{sec:thealgorithm} to construct the DtN map $T$ for the problem
$$
-\Delta u(\pxx) - \kappa^{2}u(\pxx) + \kappa^{2}b(\pxx)\,u(\pxx) = 0,\qquad \pxx \in \Omega.
$$
Then we know that
\begin{equation}
\label{eq:comstock1}
v_{n}|_{\Gamma} + w_{n}|_{\Gamma} = T\,\bigl(v|_{\Gamma} + w|_{\Gamma}\bigr),
\end{equation}
where $v_{n}$ and $w_{n}$ are normal derivatives of $v$ and $w$,
respectively.  Now use BIE methods to construct the DtN map $S$ for
the problem (\ref{eq:comstock2}) on the exterior domain $\Omega^{\rm c}$. Then $w$ must satisfy
(\ref{eq:comstock2}).
\begin{equation}
\label{eq:comstock3}
w_{n}|_{\Gamma} = S\,w|_{\Gamma}.
\end{equation}
Combining (\ref{eq:comstock1}) and (\ref{eq:comstock3}) we find
$$
v_{n}|_{\Gamma} + S\,w|_{\Gamma} = T\,v|_{\Gamma} + T\,w|_{\Gamma}.
$$
In other words,
\begin{equation}
\label{eq:comstock4}
(S - T)\,w|_{\Gamma} = T\,v|_{\Gamma} - v_{n}|_{\Gamma}.
\end{equation}
Observing that both $v|_{\Gamma}$ and $v_{n}|_{\Gamma}$ can be
obtained from (\ref{eq:v}), and that $S$ and $T$ are now available, we
see that $w|_{\Gamma}$ can be determined by solving
(\ref{eq:comstock4}).

\subsection{Problems in three dimensions}
There is no conceptual difficulty in generalizing the method to
problems in $\mathbb{R}^{3}$.  However, since the fraction of points
located on interfaces will increase, the complexity of the
pre-computation and the solve stages will be $O(N^{2})$ and
$O(N^{4/3})$, respectively.  The acceleration techniques described in
Section \ref{sec:linear} will likely become very valuable in
constructing highly efficient implementations in 3D.

\subsection{Formulating a scheme impervious to resonances}
\label{sec:reformulation}
As discussed in Section \ref{sec:resonance}, the fact that the proposed
scheme relies on Dirichlet-to-Neumann maps causes problems when the
hierarchical partitioning of the domain involves a sub-domain that admits
resonant modes. While this issue can be managed quite easily, it would
clearly be preferable to formulate a variation of the scheme that is
inherently not vulnerable in this regard. One possible approach would be
use a so called ``total wave'' approach suggested by Yu Chen of New York
University (private communication). The idea is to maintain for each leaf not an
operator that maps Dirichlet data to Neumann data, but rather a
collection of matching pairs of Dirichlet and Neumann data
(represented as vectors of tabulated values on the boundary). If you
know the collection of pairs for two adjacent boxes, you can construct
the collection for the union box via a merge procedure similar to the
one used in Section \ref{sec:merge}. This approach appears to not
suffer from any problems in the case where one of the involved boxes
admits resonant modes, but has the drawback that it would be less amenable
to the acceleration technique described in Section \ref{sec:linear}.


\section{Conclusions}
\label{sec:conc}

The paper describes a composite spectral scheme for solving variable coefficient
elliptic PDEs with smooth coefficients on simple domains such as squares and rectangles.
The method involves a \textit{direct} solver and can in a single sweep solve problems
for which state-of-the-art iterative methods require thousands of iterations.
High order spectral approximations are used. As a result, potential fields can be
computed to a relative precision of about $10^{-10}$ using twelve points per wave-length
or less.

Numerical experiments indicate that the method is very fast. For a problem involving
1.6M degrees of freedom discretizing a domain of size $100 \times 100$ wavelengths,
the pre-computation stage of the direct solver took less than 2 minutes on a laptop.
Once the solution operator had been computed, the actual solve that given a vector of
Dirichlet data on the boundary constructs the solution at all 1.6M internal tabulation
points required only 0.3 seconds. The computed solution had a relative accuracy
of $10^{-9}$.

The asymptotic complexity of the method presented is $O(N^{1.5})$ for the construction of
the solution operator, and $O(N \log N)$ for a solve once the solution operator has been
created. For a situation where $N$ is increased while the wave-number is kept fixed, it
appears possible to improve the asymptotic complexity to $O(N)$ for both the pre-computation
and the solve stages (see Section \ref{sec:linear}), but such a code has not yet been written.

The method presented has a short-coming in that it is in principle vulnerable to
resonances. It relies on a hierarchical partitioning of the domain, and if any one
of the boxes in this partitioning is resonant, the method breaks down. In practice,
this problem appears to happen very rarely, and can be both detected and remedied
if it does occur.

\lsp

\noindent
\textbf{Acknowledgments:} The author benefitted greatly from
conversations with Vladimir Rokhlin of Yale university. Valuable
suggestions were also made by Yu Chen (NYU), Adrianna Gillman
(Dartmouth), Leslie Greengard (NYU), and Mark Tygert (NYU). The work
was supported by the NSF under contracts 0748488 and 0941476, and by
the Wenner-gren foundation.  Most of the work was conducted during a
sabbatical spent at the Department of Mathematics at Chalmers
University and at the Courant Institute at NYU. The support from these
institutions is gratefully acknowledged.




\bibliography{main_bib}
\bibliographystyle{amsplain}

\begin{appendix}
\section{Graphical illustration of the hierarchical merge process}
\label{app:cartoon}

This section provides an illustrated overview of the hierarchical
merge process described in detail in Section \ref{sec:thealgorithm}
and in Figure \ref{fig:precomp}. The figures illustrate a situation
in which a square domain $\Omega = [0,1]^{2}$ is split into $4 \times 4$
leaf boxes on the finest level, and an $8\times 8$ spectral grid is
used in each leaf.

\vspace{2mm}

\noindent
\textbf{Step 1:} Partition the box $\Omega$ into $16$ small boxes that
each holds an $8\times 8$ Cartesian mesh of Chebyshev nodes.
For each box, identify the internal nodes (marked in white) and
eliminate them as described in Section \ref{sec:leaf}. Construct
the solution operator $\UU$, and the DtN operators
encoded in the matrices $\VV$ and $\WW$.

\begin{center}
\setlength{\unitlength}{1mm}
\begin{picture}(130,57)
\put(00,00){\includegraphics[height=55mm]{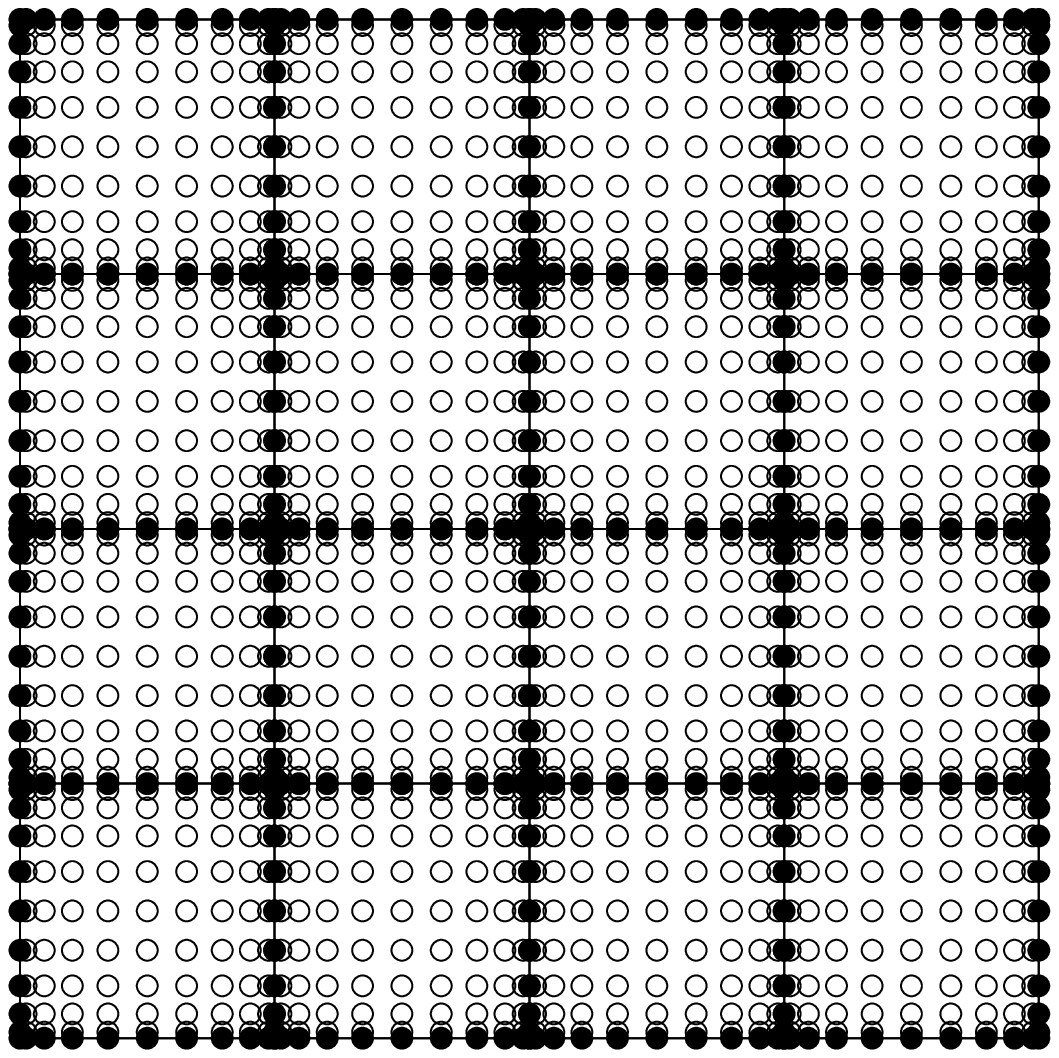}}
\put(63,27){$\Rightarrow$}
\put(60,32){Step 1}
\put(73,00){\includegraphics[height=55mm]{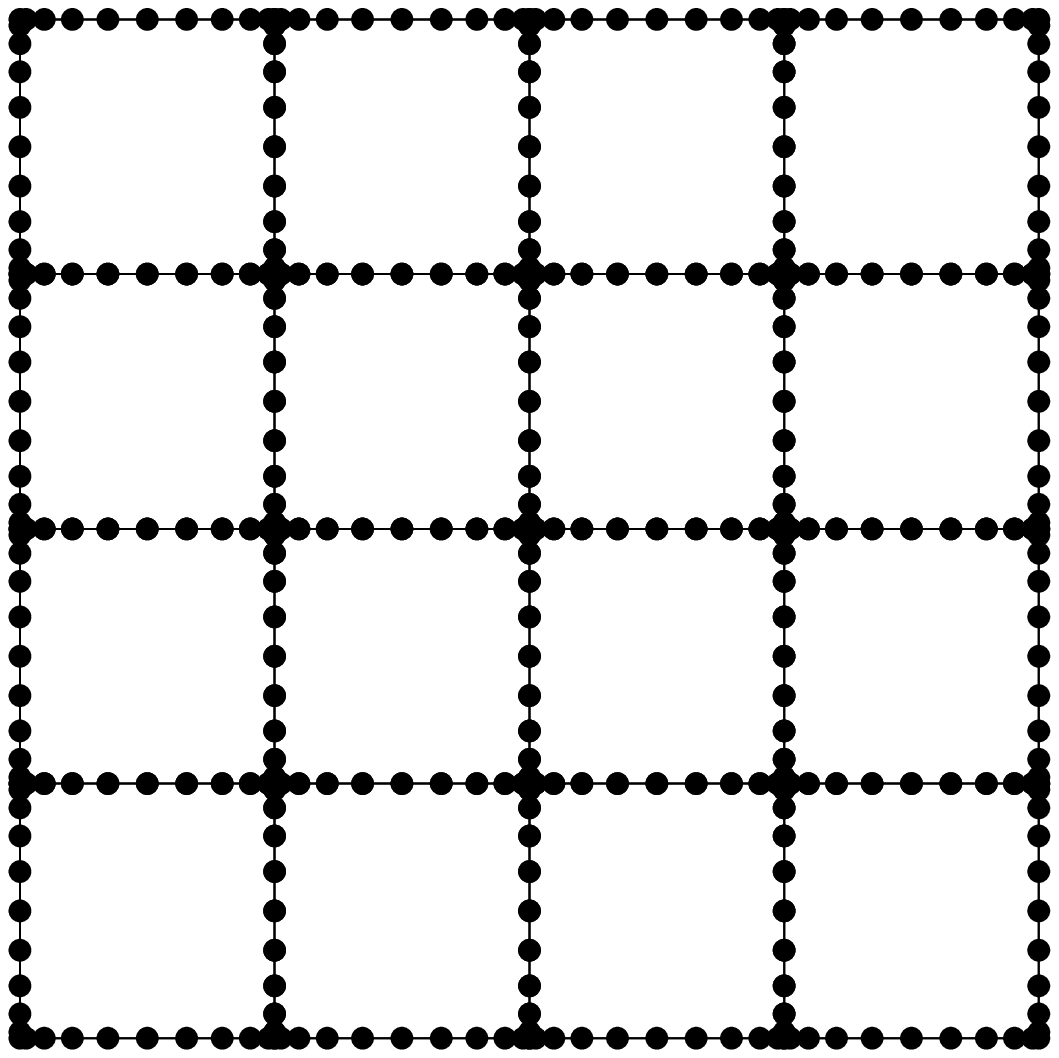}}
\end{picture}
\end{center}

\vspace{2mm}

\noindent
\textbf{Step 2:} Merge the small boxes by pairs as described in Section \ref{sec:merge}.
The equilibrium equation for each rectangle is formed using the DtN operators of the two
small squares it is made up of. The result is to eliminate the interior nodes (marked in
white) of the newly formed larger boxes. Construct the solution operator $\UU$ and
the DtN matrices $\VV$ and $\WW$ for the new boxes.

\begin{center}
\setlength{\unitlength}{1mm}
\begin{picture}(130,57)
\put(00,00){\includegraphics[height=55mm]{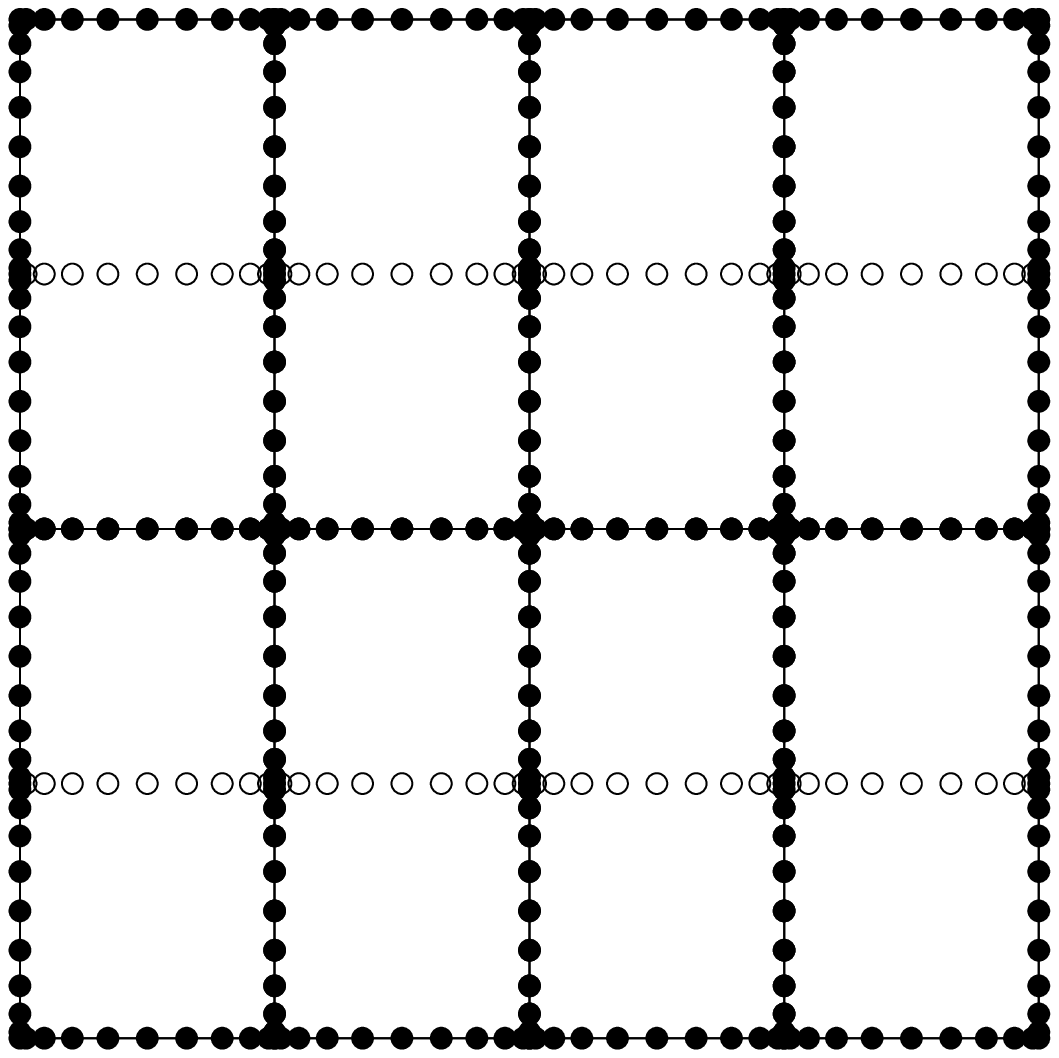}}
\put(63,27){$\Rightarrow$}
\put(60,32){Step 2}
\put(73,00){\includegraphics[height=55mm]{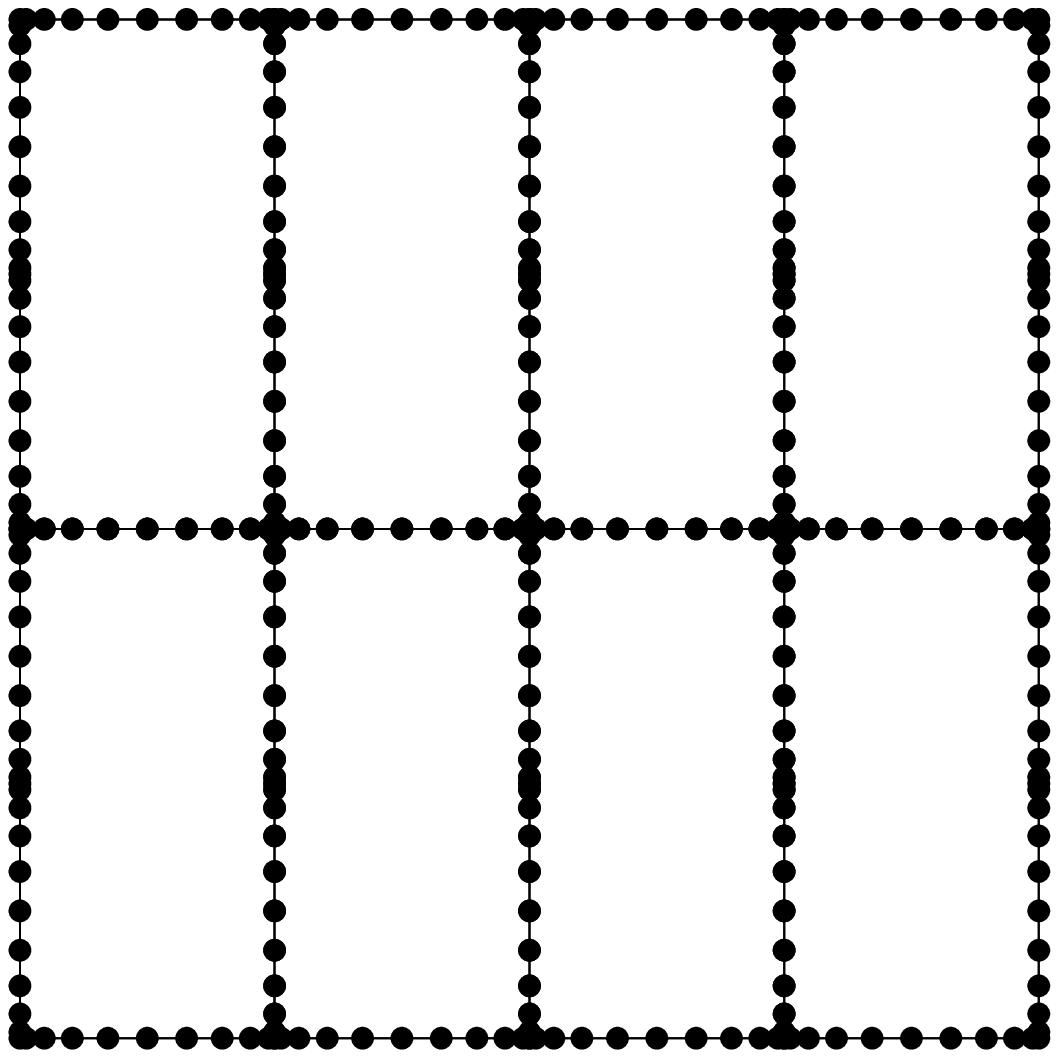}}
\end{picture}
\end{center}

\vspace{2mm}

\noindent
\textbf{Step 3:} Merge the boxes created in Step 2 in pairs, again via
the process described in Section \ref{sec:merge}.

\begin{center}
\setlength{\unitlength}{1mm}
\begin{picture}(130,57)
\put(00,00){\includegraphics[height=55mm]{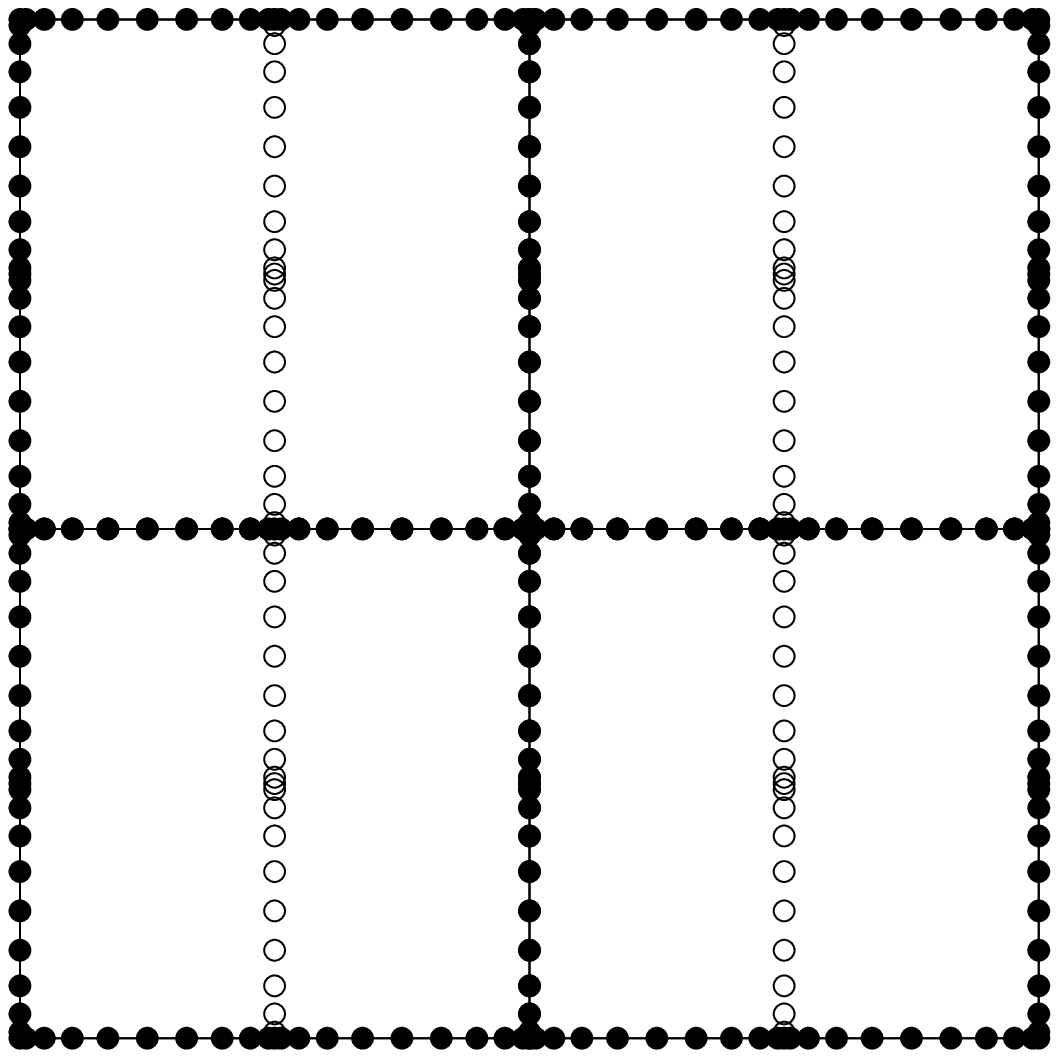}}
\put(63,27){$\Rightarrow$}
\put(60,32){Step 3}
\put(73,00){\includegraphics[height=55mm]{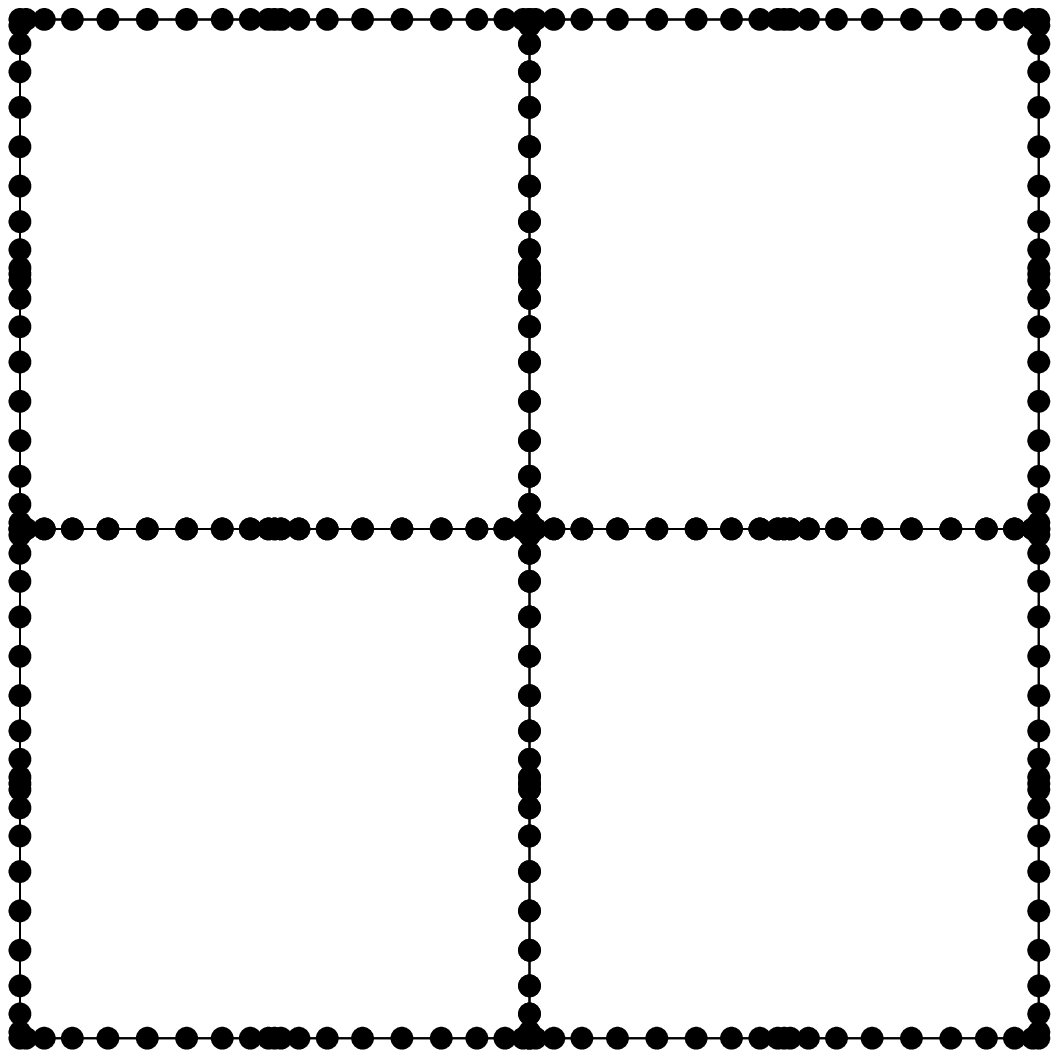}}
\end{picture}
\end{center}

\vspace{2mm}

\noindent
\textbf{Step 4:} Repeat the merge process once more.

\begin{center}
\setlength{\unitlength}{1mm}
\begin{picture}(130,57)
\put(00,00){\includegraphics[height=55mm]{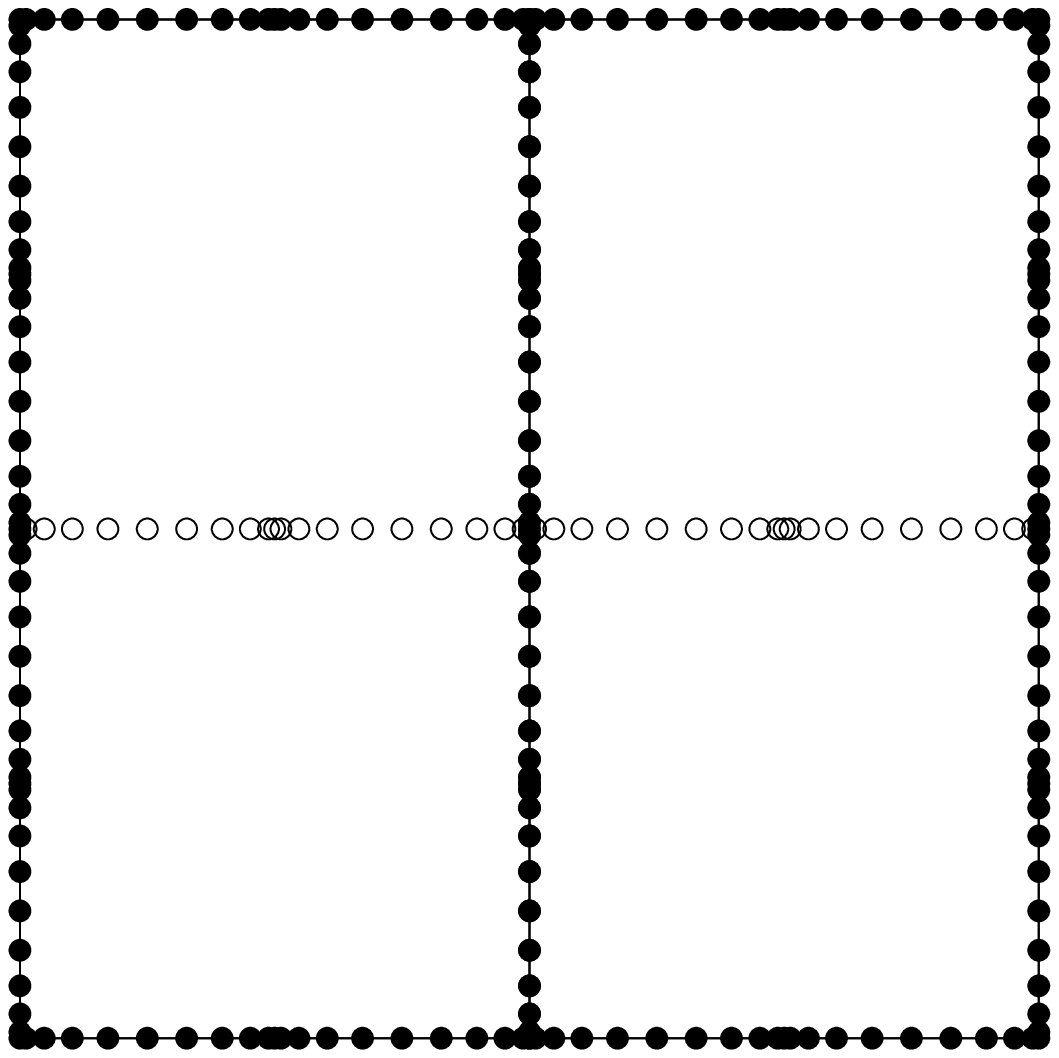}}
\put(63,27){$\Rightarrow$}
\put(60,32){Step 4}
\put(73,00){\includegraphics[height=55mm]{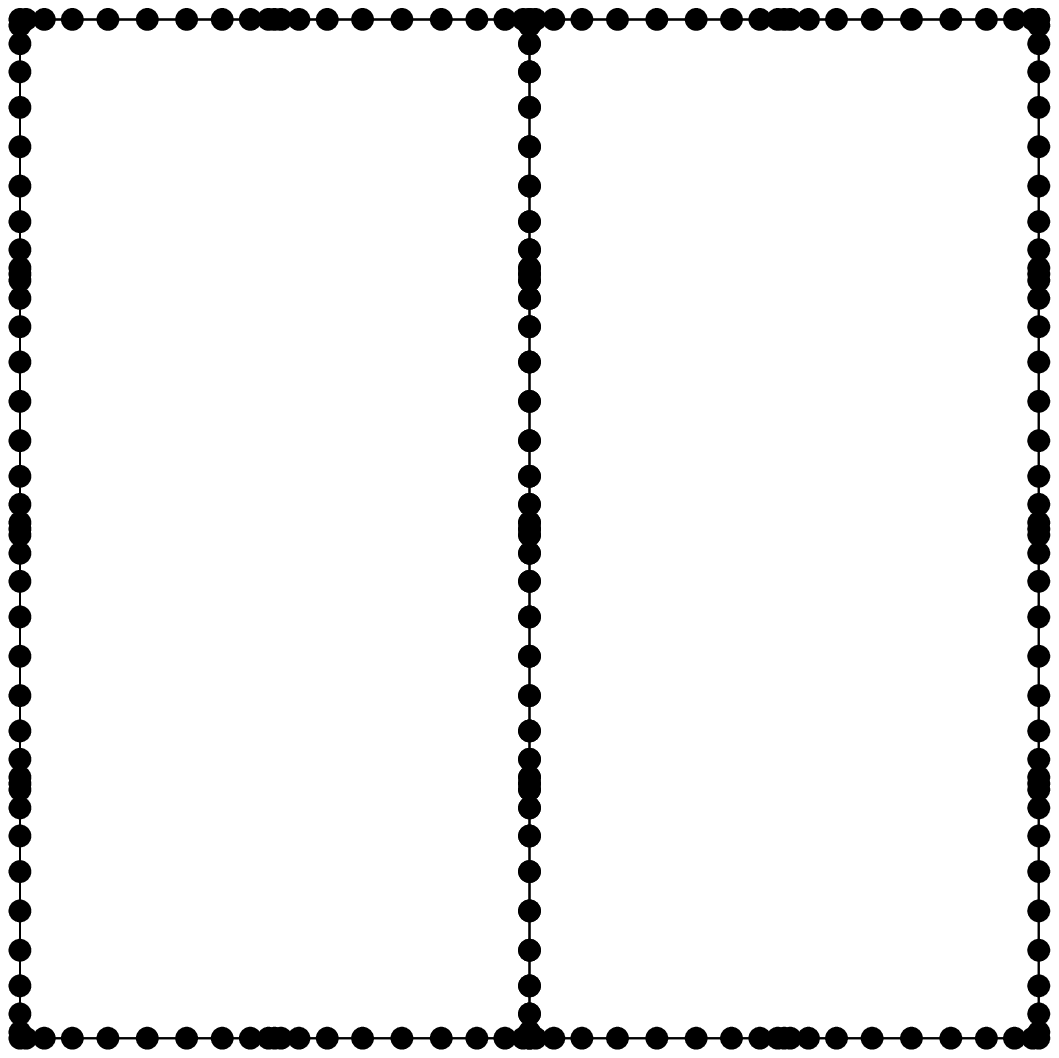}}
\end{picture}
\end{center}

\vspace{2mm}

\noindent
\textbf{Step 5:} Repeat the merge process one final time to obtain the
DtN operator for the boundary of the whole domain.

\begin{center}
\setlength{\unitlength}{1mm}
\begin{picture}(130,57)
\put(00,00){\includegraphics[height=55mm]{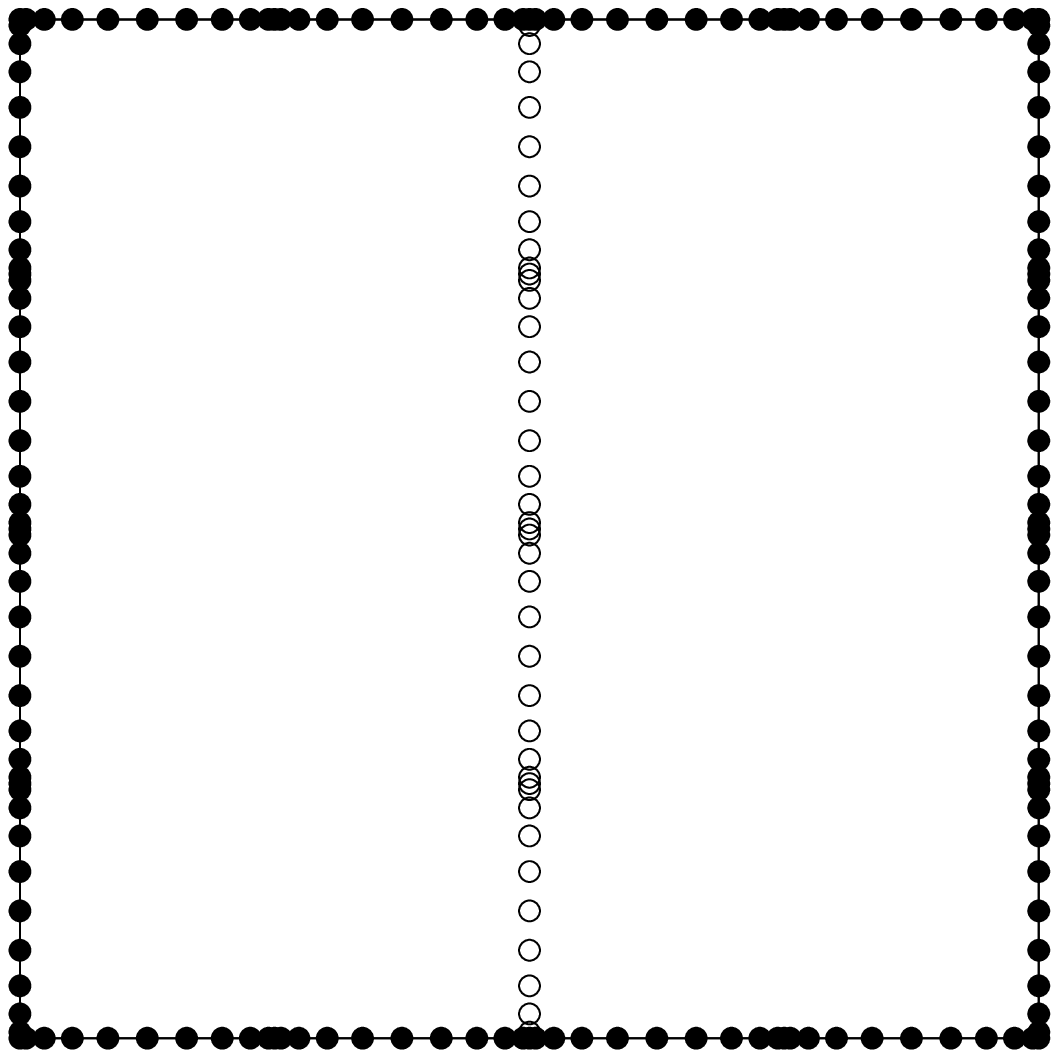}}
\put(63,27){$\Rightarrow$}
\put(60,32){Step 5}
\put(73,00){\includegraphics[height=55mm]{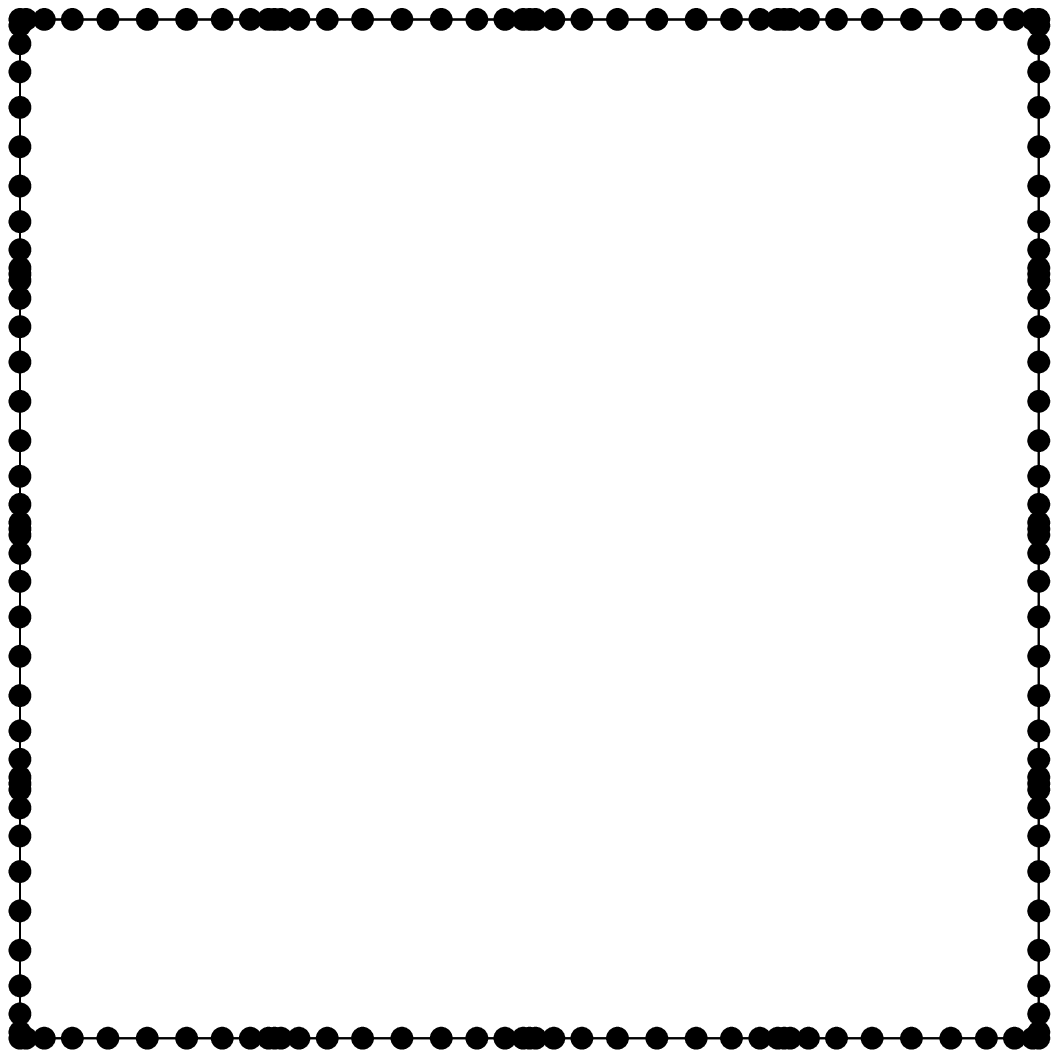}}
\end{picture}
\end{center}

\end{appendix}

\end{document}